\documentclass[12pt,a4paper]{article}

\usepackage{amsmath,amsthm,amssymb}
\usepackage{mathrsfs}
\usepackage{natbib}
\usepackage[dvips]{graphicx}

\newtheorem{thm}{Theorem}
\newtheorem{lem}{Lemma}

\newtheorem{cor}{Corollary}
\newtheorem{definition}{Definition}
\newtheorem{remark}{Remark}
\newtheorem{example}{Example}

\newtheorem{exercise}{Exercise}
\newtheorem{report}{Report}
\newtheorem{assumption}{Assumption}

\newcommand{\Proof}{\noindent {\bf Proof: }\ }

\newcommand{\bin}{\mathrm{Bin}}

\newcommand{\prob}{\mathrm{Prob}}
\newcommand{\abs}[1]{\vert #1 \vert}
\newcommand{\abslr}[1]{\left\vert #1 \right\vert}
\newcommand{\norm}[1]{\| #1 \|}
\newcommand{\normlr}[1]{\left\| #1 \right\|}
\newcommand{\real}{\mathbb{R}}
\newcommand{\dist}{\mathrm{dist}}
\newcommand{\indicator}{\mathrm{1}}
\newcommand{\rmd}{\mathrm{d}}
\newcommand{\rmor}{\mathrm{or}}

\newcommand{\locparam}{\mu}
\newcommand{\scparam}{\sigma}
\newcommand{\wtparam}{\alpha}

\newcommand{\seqnum}{j}
\newcommand{\scra}{\mathscr{A}}
\newcommand{\scrb}{\mathscr{B}}

\newcommand{\scrg}{\mathscr{G}}
\newcommand{\scrk}{\mathscr{K}}
\newcommand{\scrl}{\mathscr{L}}
\newcommand{\scrkrest}{\scrk_{{R}}}
\newcommand{\scrkbbound}{\scrk_{c_{0}<\scparam <B}}
\newcommand{\ndownarrow}{>}

\title{Strong consistency of MLE for
finite mixtures of location-scale distributions 
when the scale parameters are exponentially small}
\author{Kentaro Tanaka\thanks{Department of Industrial Engineering and Management, 
Tokyo Institute of Technology, 
2-12-1 O-okayama, Meguro-ku, Tokyo 152-8550 JAPAN, 
E-mail:\,tanaken@me.titech.ac.jp} 
\, and\, 
Akimichi Takemura\thanks{Department of Mathematical Informatics,  
Graduate School of Information Science and Technology,  
University of Tokyo, 
7-3-1 Hongo, Bunkyo-ku, Tokyo 113-8656, JAPAN, 
E-mail:\,takemura@stat.t.u-tokyo.ac.jp}
}
\date{}

\begin{document}
\maketitle \vspace{-0.65cm}
\begin{abstract}
  In a finite mixture of location-scale distributions maximum likelihood
  estimator does not exist because of the unboundedness of the
  likelihood function when the scale parameter of some mixture component
  approaches zero.  In order to study the strong consistency of
  maximum likelihood estimator, we consider the case that the scale
  parameters of the component distributions are restricted from below
  by $c_n$, where $\{c_n\}$ is a sequence of positive real numbers which tend
  to zero as the sample size $n$ increases.  We prove that under mild
  regularity conditions maximum likelihood estimator is
  strongly consistent if the scale parameters are restricted from below by
  $c_{n} = \exp(-n^d)$, $0 < d < 1$.
\end{abstract}

\section{Introduction}

In some finite mixture distributions maximum likelihood estimator (MLE)
does not exist.  Let us consider the following example.  Denote a
normal mixture distribution with $M$ components and parameter
$\theta=(\wtparam_1, \mu_{1}, \sigma_{1}^{2} ,\ldots,\wtparam_{M},
\mu_{M}, \sigma_{M}^{2})$ by
\begin{eqnarray}
 f(x;\theta) 
  = 
  \sum_{m=1}^{M}
  \wtparam_{m} \phi_{m}(x; \mu_{m}, \sigma_{m}^{2}) , 
 \nonumber 
\end{eqnarray}
where $\wtparam_{m} \; (m=1,\ldots,M)$ 
are nonnegative real numbers that sum to one and 
$\phi_{m}(x; \mu_{m}, \sigma_{m}^2)$ 
are normal densities.
Let $x_{1},\ldots,x_{n}$ denote 
a random sample of size $n \geq 2$ 
from the density $f(x;\theta_{0})$. 
In view of the identifiability problem of mixture distributions
discussed below, here $\theta_{0}$ is a parameter value designating the
true distribution.  However for simplicity we just say 
$\theta_{0}$ is the true parameter from now on.
The log likelihood function is 
\begin{eqnarray}
 \sum_{i=1}^{n}
  \log{f(x_{i};\theta)}
  = 
  \sum_{i=1}^{n}
  \log{
  \left\{
   \sum_{m=1}^{M}
   \wtparam_{m} \phi_{m}(x_{i}; \mu_{m}, \sigma_{m}^{2})
  \right\}
  }.
 \nonumber 
\end{eqnarray}
If we set $\mu_{1} = x_{1}$, 
then the likelihood tends to infinity as  
$\sigma_{1}^{2} \rightarrow 0$.
Thus MLE 
does not exist.

But when we restrict 
$\sigma_{m} \geq c \;(m=1,\ldots,M)$
by some positive real constant $c$, 
we can avoid the divergence of the likelihood.
Furthermore, 
in this situation, 
it can be shown that MLE 
is strongly consistent 
if the true parameter $\theta_0$ 
is in the restricted parameter space.

On the other hand, 
the smaller $\sigma_{1}^{2}$ is,
the less contribution $\phi_{1}(x; \mu_{1}=x_1, \sigma_{1}^{2})$
makes to 
the likelihood at other observations $x_{2},\ldots,x_{n}$.
Therefore an interesting question here is whether we 
can decrease the bound $c=c_n$ to zero with the sample size $n$ and
yet guarantee the strong consistency of MLE.
If this is possible, the further question is how fast
$c_n$ can decrease to zero. 

This question is similar to the (so far open) 
problem stated in \citet{H1985}, 
which treats mixtures of normal distributions 
with constraints imposed on the ratios of variances 
while our restriction is imposed on variances themselves.
See also a discussion in section 3.8.1 of \citet{M2000}.

In the above example, the normality of the component distributions is
not essential and the same difficulty exists for finite mixtures of
general location-scale distributions such as mixtures of uniform
distributions. Furthermore in this paper we allow that each component
belongs to different location-scale families.  Let $\scparam_{m}$
($m=1,\ldots,M$) denote the scale parameters of the component
distributions and consider the restriction $\scparam_{m} \geq c_{n} \;
(m = 1, \ldots, M)$.  Then a question of interest here is whether we
can decrease the bound $c_n$ to zero.

For the case of mixture of uniform distributions, in \citet{TT2003-20} 
we proved that MLE 
is strongly consistent if $c_{n} = \exp(-n^d)$, $0< d<1$.  Here $d$
can be arbitrarily close to 1 but fixed.
In this paper, we prove that the same result holds for general finite
mixtures 
of location-scale distributions under very mild regularity conditions
(assumptions  \ref{assumption:1}--\ref{assumption:4} below).
We employ the same line of proof as in 
\citet{TT2003-20}, 
but the proof for the general finite
mixture is much more difficult.
As discussed in section
\ref{sec:discussion} the normal density satisfies the
regularity conditions and our result implies that 
MLE is strongly
consistent for the finite normal mixture
if $\sigma_m \ge c_n = \exp(-n^d)$, $0<d<1$, $m=1,\ldots,M$.

Our framework is closely related to the method of sieve
(\citet{G1981}).
In the sieve method 
an objective function  is maximized 
over a constrained subspace of parameter space 
and then this subspace is expanded to the whole parameter space 
as the sample size increases.
Some applications and 
consistency results for the method 
are given in \citet{G1982}.
MLE based on a sieve is called a sieve MLE. 
The convergence rates of sieve MLE 
for Gaussian mixture problems are studied in 
\citet{GW2000} and \citet{GV2001} 
and their ideas are very interesting.
They obtain the convergence rates by 
bounding the Hellinger bracketing entropy of 
subsets of the function space 
and assume that the corresponding subsets of the parameter space are compact 
so that their bracketing entropy does not diverge.
In the case of sieve MLE, the approximating subspaces are usually
taken to be compact, whereas 
we treat a sequence of non-compact subsets of the parameter space 
expanding  to the whole parameter space as the sample size increases.
Therefore results on sieve MLE are not directly applicable in
our framework. 

The organization of the paper is as follows.
In section
\ref{sec:preliminaries} we summarize some preliminary descriptions.
In section \ref{sec:main} we state our main results in theorems
\ref{kappa-lambda:thm} and \ref{main-thm:thm}.
Section \ref{sec:proof}
is devoted to the proof of theorems and lemmas. 
Finally in section \ref{sec:discussion} we
give some discussions.

\section{Preliminaries on strong consistency and identifiability of mixture distributions}
\label{sec:preliminaries}

A mixture of $M$ densities with parameter 
$\theta=(\wtparam_1, \locparam_{1}, \scparam_{1},\ldots,\wtparam_{M}, \locparam_{M}, \scparam_{M})$ 
is defined by 
\begin{eqnarray}
 f(x;\theta) \equiv \sum_{m=1}^{M} \wtparam_{m} f_{m}(x;\locparam_{m}, \scparam_{m}) , 
  \nonumber 
\end{eqnarray}
where $\wtparam_{m}$, $m=1,\ldots,M$, called the mixing weights, 
are nonnegative real numbers that sum to one 
and $f_{m}(x;\locparam_{m}, \scparam_{m})$, called the components of the mixture,
are density functions.  In this paper we consider the case that 
the component densities are location-scale densities 
with the location parameter $\locparam_{m}\in \mathbb{R}$
and the scale parameter $\scparam_{m}>0$, i.e. 
\begin{equation}
\label{eq:location-scale}
f_{m}(x;\locparam_{m}, \scparam_{m})=\frac{1}{\scparam_m} f_m \left( \frac{x-\locparam_m}{\scparam_m}; 0,1\right).
\end{equation}
As mentioned above, we allow $f_{m}(x;\locparam_{m}, \scparam_{m})$ to
belong to different families.
For example, $f_{1}(x;\locparam_{1}, \scparam_{1})$ may be  a normal density, 
$f_{2}(x;\locparam_{2}, \scparam_{2})$ may be a uniform density, etc.
Let $\Omega_{m}=\mathbb{R}\times (0,\infty)$ denote 
the parameter space of the  $m$-th component $(\locparam_m, \scparam_m)$ 
and let $\Theta$ denote the entire parameter space: 
\begin{eqnarray}
 \Theta  & \equiv & 
 \{
 (\wtparam_{1}, \ldots , \wtparam_{M}) \in \real^{M}
 \mid 
 \sum_{m=1}^{M} \wtparam_{m} = 1 
 \; , \;
 \wtparam_{m} \geq 0
 \}\times \prod_{m=1}^{M}\Omega_{m}.
 \nonumber 
\end{eqnarray}

Let $\scrk$ be a subset of $\{1,2,\ldots,M\}$ and 
let $\abslr{\scrk}$ denote the number of elements in
$\scrk$. 
Denote by $\theta_{\scrk}$ a subvector of $\theta \in \Theta$ 
consisting of the components in $\scrk$. 
Then the parameter space of subprobability measures 
consisting of the components in $\scrk$ 
is
\begin{equation}
 \bar{\Theta}_{\scrk} \equiv 
  \{
  \theta_{\scrk} \mid \theta \in \Theta, \sum_{m \in \scrk} \alpha_m
  \le 1 
  \}.
  \label{eq:def:barTheta_scrk}
\end{equation}
Corresponding density and the set of subprobability densities are
denoted by 
\begin{align}
 f_{\scrk}(x;\theta_{\scrk})
  &\equiv \sum_{k \in \scrk} \wtparam_{k} f_k(x;{\locparam}_k, \scparam_{k}) ,
  \label{eq:def:f_scrk} \\
 \scrg_{\scrk} 
  & \equiv 
  \{
  f_{\scrk}(x;\theta_{\scrk}) 
  \mid 
  \theta_{\scrk} \in \bar{\Theta}_{\scrk}
  \} .
  \label{eq:def:scrg_scrk}
\end{align}
Furthermore denote 
the set of subprobability densities with no more than $K$
components by
\begin{equation}
\scrg_{K} \equiv \bigcup_{\abslr{\scrk}\leq K}\scrg_{\scrk}
\qquad (1 \leq K \leq M). 
  \label{eq:def:scrg_K}
\end{equation}

We now briefly discuss identifiability of parameters.
In mixture models, different parameters may designate the same
distribution.  When the component densities belong to a common
location-scale family, we can permute the labels of the components
and the distribution remains the same.  A mixture model of $K-1$
components can be obtained by setting one weight $\alpha_m = 0$ (with
arbitrary $\mu_m$ and $\sigma_m$) in a
model with  $K$ components.  These are trivial cases of unidentifiability
of parameters.
However there are more
complicated cases.  Let $U(x;a,b)$ denote the uniform density on the
interval $[a,b]$. Then, for example,
$\frac{1}{3}U(x;-1,1) + \frac{2}{3}U(x;-2, 2)$
and 
$\frac{1}{2}U(x;-2,1) + \frac{1}{2} U(x;-1,2)$
represent the same distribution (\cite{EH1981}). 
In this case the limiting behavior of MLE is not obvious,
although the estimated density should be consistent. 
Therefore we first give a definition of consistency in terms of 
the estimated density.

Let $f_0(x) = f(x;\theta_0)$ denote the true density and let
$\hat f_n(x)=f(x;\hat\theta_n)$ denote the estimated density.  
\begin{definition}
\label{definition:1}
An estimator $\hat{f}_{n}$ 
is strongly consistent if
\begin{equation}
 \prob \left(
  \lim_{n \rightarrow \infty}
  \normlr{\hat{f}_{n}-f_{0}} = 0
 \right) = 1 \; , 
\label{strongconsistency:def}
\nonumber 
\end{equation}
where 
$\norm{\cdot}$ is the $L_{1}$-norm. 
\end{definition}

Although definition \ref{definition:1} is conceptually simple, in
order to prove the strong consistency of MLE
we work with the location and the scale parameters in
(\ref{eq:location-scale}) and the mixing weights.
In order to deal with the identifiability problem 
let us introduce a distance between 
two sets of parameters. Let 
$\dist(\theta,\theta')$ denote the ordinary Euclidean distance (or any
other equivalent distance) between
two parameter vectors $\theta,\theta' \in \Theta$. 
For $U,V \subset \Theta$ define
\begin{eqnarray}
 \dist(U,V) 
  \equiv 
  \inf_{\theta \in U} \inf_{\theta' \in V} \dist(\theta,\theta') .
  \nonumber 
\end{eqnarray}
For a parameter $\theta$, let 
\begin{eqnarray}
 \Theta(\theta) 
  \equiv 
  \{\theta' \in \Theta \mid f(x;\theta') = f(x;\theta)
  \quad 
  \forall 
  x\} . 
  \nonumber
\end{eqnarray}
Then $\Theta_{0} = \Theta(\theta_{0})$ denotes the set of true parameters. 
Since our densities are continuous with respect to $\theta$, 
by Scheff\'e's theorem~(Theorem 16.12  of \cite{B1995}) 
$\dist(\Theta(\hat\theta_n), \Theta_0) \rightarrow 0$  implies
$\normlr{\hat{f}_{n}-f_{0}} \rightarrow 0$.

\section{Main results}
\label{sec:main}

We assume the following regularity conditions for 
strong consistency of MLE.
\begin{assumption}
 There exist real constants 
 $v_{0},v_{1} > 0$ and $\beta > 1$ such that 
 \begin{eqnarray}
  f_{m}(x;\locparam_{m}=0, \scparam_{m}=1) 
   \leq 
   \min \{v_{0} \;, \; v_{1} \cdot \abs{x}^{-\beta}\}
   \nonumber 
 \end{eqnarray}
 for all $m$.
\label{assumption:1}
\end{assumption}

This assumption means that $f_{m} \; (m=1,\ldots, M)$ are bounded and
their tails decrease to zero faster than or equal to
$\abs{x}^{-\beta}$, which is a very mild condition.

The following three regularity conditions are standard 
conditions assumed in discussing strong consistency of MLE.
Let $\Gamma$ denote any compact subset of $\Theta$.

\begin{assumption}
 For $\theta \in \Theta$ and any positive real number $\rho$, let 
 \begin{eqnarray}
  f(x;\theta,\rho) 
   & \equiv & \sup_{\dist(\theta',\theta) \leq \rho}f(x;\theta')
   . 
   \nonumber
 \end{eqnarray}
 For each $\theta \in \Gamma$ and sufficiently small $\rho$, 
 $f(x;\theta, \rho)$ is measurable. 
\label{assumption:2}
\end{assumption}
\begin{assumption}
 For each $\theta \in \Gamma$,  if\/
 $\lim_{\seqnum \rightarrow \infty}\theta^{(\seqnum)} = \theta, \; (\theta^{(\seqnum)} \in \Gamma)$ 
 then\/
 $\lim_{\seqnum \rightarrow \infty}
 f(x;\theta^{(\seqnum)}) = f(x;\theta)$
 except on a set which is a null set and 
 does not depend on the sequence 
 $\{\theta^{(\seqnum)}\}_{\seqnum=1}^{\infty}$.
\label{assumption:3}
\end{assumption}
\begin{assumption}
 \begin{eqnarray}
  \int \abslr{
   \log{f(x;\theta_0)}
   }
   f(x;\theta_0)\rmd x < \infty .
   \nonumber 
 \end{eqnarray}
\label{assumption:4}
\end{assumption}

Let 
$E_{0}[\cdot]$ denote the expectation under the true parameter $\theta_{0}$.
The following theorem is essential to our 
argument and it is of some independent interest.
\begin{thm}
 \label{kappa-lambda:thm}
 Suppose that assumptions \ref{assumption:1}--\ref{assumption:4}
 are satisfied 
 and $f_{0} \in \scrg_{M} \backslash \scrg_{M-1}$ 
 where $\scrg_{M}$ and $\scrg_{M-1}$ are defined in
$(\ref{eq:def:scrg_K})$. 
Then there  exist real constants $\kappa, \lambda > 0$ such that  
 \begin{eqnarray}
  E_{0}
  \left[
   \log{
   \left\{ 
    g(x) + \kappa
   \right\}
   }
  \right]
  +
  \lambda
  <
  E_{0}[\log{f(x;\theta_{0})}] 
  \label{kappa-lambda:eq}
 \end{eqnarray}
 for all 
 $g \in \scrg_{M-1}$. 
\end{thm} 

We now state the main theorem of this paper.
\begin{thm}
 Suppose that assumptions \ref{assumption:1}--\ref{assumption:4}
 are satisfied 
 and $f_{0} \in \scrg_{M}\backslash \scrg_{M-1}$
 where $\scrg_{M}$ and $\scrg_{M-1}$ are defined in
$(\ref{eq:def:scrg_K})$. 
 Let $c_{0} > 0$ and $ 0 < d < 1$.
 If $c_{n} = c_{0} \cdot \exp(-n^{d})$ and 
 \begin{eqnarray}
  \Theta_{n} 
   \equiv 
   \{
   \theta \in \Theta
   \mid 
   \scparam_{m} \geq c_{n}\; , \; (m=1,\ldots,M)
   \},
   \nonumber 
 \end{eqnarray}
 then
  \begin{eqnarray}
  \prob
   \left(
    \lim_{n \rightarrow \infty}
    \dist(\Theta(\hat{\theta}_{n}), \Theta_0) = 0
   \right) = 1 \; ,  
   \nonumber 
\end{eqnarray}
where $\hat{\theta}_{n}$ is MLE
restricted to $\Theta_{n}$.  
 \label{main-thm:thm}
\end{thm}

As remarked at the end of the previous section
theorem \ref{main-thm:thm} implies 
the following corollary. 
\begin{cor}
Under the same assumptions of $theorem~\ref{main-thm:thm}$, 
$\hat{f}_{n}$ is strongly consistent 
in the sense of definition~\ref{definition:1}. 
\end{cor}

\section{Proofs}
\label{sec:proof}

In this section, we prove theorems stated in section \ref{sec:main}.
The organization of this section is as follows.  First in
subsection~\ref{sec:notat-some-lemm} we state some lemmas for
theorem~\ref{kappa-lambda:thm} and \ref{main-thm:thm}.  Next 
in subsection~\ref{sec:proof-kappa-lambda-thm}
we prove
theorem~\ref{kappa-lambda:thm} 
which is also essential for theorem~\ref{main-thm:thm}.
Finally we prove theorem~\ref{main-thm:thm} 
in subsection~\ref{sec:proof-main-thm}. 
For convenience a list of
notations used in our proofs is provided at the end of this paper.

\subsection{Notations and some lemmas}
\label{sec:notat-some-lemm}

\begin{figure}[htbp]
 \begin{center}
  \includegraphics{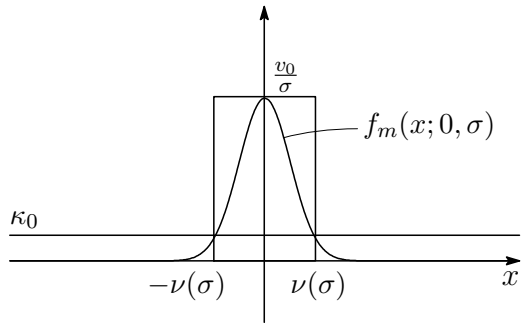}
  \caption{Each component is bounded by step function.}
  \label{fig:BoundedByStepFunction}
 \end{center}
\end{figure}

Fix arbitrary $\kappa_0> 0$, which corresponds to $\kappa$ in theorem
\ref{kappa-lambda:thm}.
Define
$\tilde{\beta}$ and  $\nu(y)$, $y>0$,  as 
\begin{eqnarray}
 \tilde{\beta}
  \equiv 
  \frac{\beta-1}{\beta}
  \quad , \quad 
 \nu(y) \equiv 
  \left(
   \frac{v_{1}}{\kappa_{0}}
  \right)^{\frac{1}{\beta}}
  y^{\tilde{\beta}}
\ ,
 \label{eq:def:beta_and_nu_and_v}
\end{eqnarray}
where 
$v_1$ and $\beta$ are given in assumption~\ref{assumption:1}.
Noting that $v_1 \cdot (\nu(y))^{-\beta} = \kappa_0 /y$, 
the following lemma is easily proved and we omit its proof. 
See figure \ref{fig:BoundedByStepFunction}.

\begin{lem}
\label{lem:BoundComponentByStepFunc}
Under the assumption~\ref{assumption:1}, 
for arbitrary $\kappa_0> 0$
each component density 
$f_{m}(x; \locparam, \scparam)$ is bounded by a step function 
\begin{eqnarray}
 f_{m}(x; \locparam, \scparam) 
  \leq 
  \max\{\indicator_{[\locparam-\nu(\scparam), \locparam+\nu(\scparam))}(x) \cdot \frac{v_0}{\scparam}
  \; , \; {\kappa_{0}}\}
  \leq 
  \indicator_{[\locparam-\nu(\scparam), \locparam+\nu(\scparam))}(x)
  \cdot 
  \frac{v_0}{\scparam}
  +
  {\kappa_{0}}, 
  \nonumber 
\end{eqnarray}
 where $\indicator_{U}(x)$ denotes the indicator function of 
$U \subset \real$.
\end{lem}

From lemma~\ref{lem:BoundComponentByStepFunc}
\begin{eqnarray}
  \sum_{m=1}^{M}
 \wtparam_{m} f_{m}(x; \locparam_{m}, \scparam_{m})
 & \leq & 
  \sum_{m=1}^{M}
  \indicator_{[\locparam_{m}-\nu(\scparam_{m}),
    \locparam_{m}+\nu(\scparam_{m}))}(x)
    \cdot 
\frac{v_0}{\scparam_{m}}
  +
  {\kappa_{0}}.
  \label{bounded-by-uniform-density:eq}
\end{eqnarray}
The right-hand side of (\ref{bounded-by-uniform-density:eq}) is a step function. 
We look at this step function where the density $f(x;\theta)$ is high,
i.e.\ the scale parameter of some component is small.

\begin{figure}[htbp]
 \begin{center}
  \includegraphics{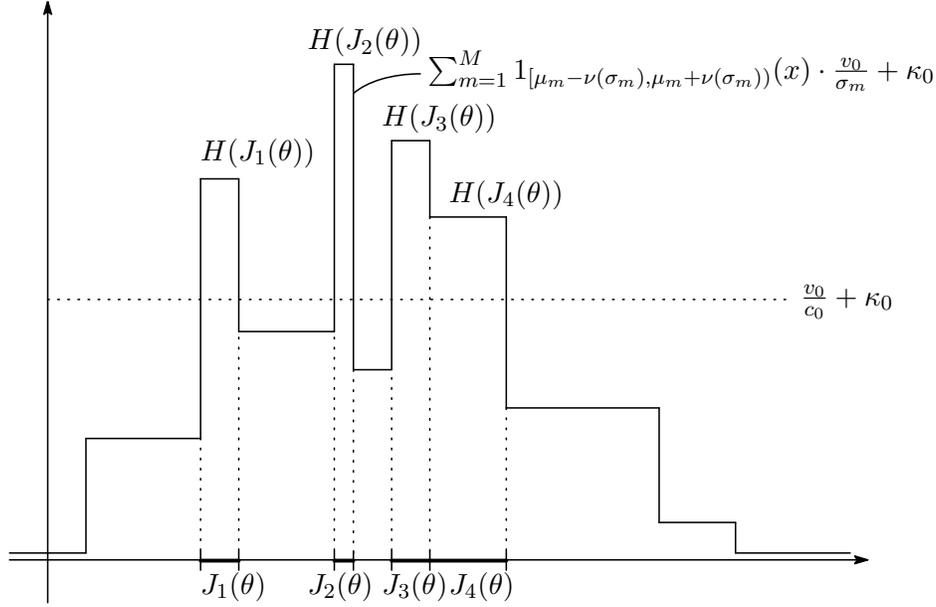}
  \caption{Definition of $J_{t}(\theta)$.}
  \label{fig:J_t}
 \end{center}
\end{figure}

For a given choice of $\kappa_0>0$, choose $c_0 > 0$ such that
\begin{equation}
 {c_0} < \frac{v_{0}}{\kappa_0(M+1)} \; , 
  \label{eq:kappa-c0-global-condition}
\end{equation}
Below  we will impose 
additional constraints on $\kappa_{0}$ and $c_{0}$ to make
$\kappa_{0}$ and $c_{0}$ sufficiently small to satisfy other conditions.
For each $\theta$, let 
\[
 \scrk_{{\scparam} \leq c_{0}} = \scrk_{{\scparam} \leq c_{0}}(\theta)
  \equiv
  \{m \mid 1 \leq m \leq M \; , \; \scparam_{m} \leq c_{0}\}  
  \nonumber 
\]
denote the set of components with $\sigma_m \le c_0$ and define
\begin{equation}
 J(\theta)
  \equiv
 \bigcup_{m \in \scrk_{{\scparam} \leq c_{0}}}
 [\locparam_{m} - \nu(\scparam_{m}), \locparam_{m} + \nu(\scparam_{m})). 
 \label{eq:def:Jtheta}
\end{equation}
On $J(\theta)$ the density $f(x;\theta)$ is high.  Now dividing
$J(\theta)$ according to the height of the step function on the right-hand side of
(\ref{bounded-by-uniform-density:eq}), for $x \in J(\theta)$ 
we can write the right-hand side of
(\ref{bounded-by-uniform-density:eq}) as
\begin{equation}
 \indicator_{J(\theta)}(x) \cdot 
\left\{
 \sum_{m=1}^{M}
  \indicator_{[{\locparam}_{m}-\nu({\scparam}_{m}), {\locparam}_{m}
    +\nu({\scparam}_{m}))}(x) \cdot 
   \frac{v_0}{{\scparam}_{m}}
  +
  {\kappa_{0}}
  \right\}
  = 
  \sum_{t=1}^{T(\theta)}
  H(J_{t}(\theta)) \cdot \indicator_{J_{t}(\theta)}(x), 
  \nonumber 
\end{equation}
where $J_{t}(\theta) \; (t=1,\ldots,T(\theta))$ are disjoint
intervals, $[\locparam_{m}-\nu(\scparam_{m}),
\locparam_{m}+\nu(\scparam_{m}))\; \; (m \in \scrk_{{\scparam} \leq
  c_{0}})$ are unions of some of $J_{t}(\theta)$'s and the height
$H(J_{t}(\theta))$ for each $t$ is defined by any $x \in
J_{t}(\theta)$ as
\begin{eqnarray}
 H(J_{t}(\theta)) \equiv  
  \sum_{m=1}^{M} \indicator_{[\locparam_{m}-\nu(\scparam_{m}),
    \locparam_{m}+\nu(\scparam_{m}))}(x) \cdot
  \frac{v_0}{\scparam_{m}}
  + {\kappa_{0}}.
\label{eq:def:HJtn}
\end{eqnarray}
See figure \ref{fig:J_t}.  For $x\in J_t(\theta)$, there is at least
one $m=m_t$ such that 
$x\in [\locparam_{m}-\nu(\scparam_{m}),
\locparam_{m}+\nu(\scparam_{m}))$ and 
$H(J_{t}(\theta))\ge v_0/c_0+\kappa_0$.
Also note that the total number $T(\theta)$ of $J_{t}(\theta)$'s satisfies
$T(\theta) \le 2M$, because the change of the height can only occur at 
$\locparam_{m}-\nu(\scparam_{m})$ or $\locparam_{m}+\nu(\scparam_{m})$.

By   (\ref{bounded-by-uniform-density:eq}) 
we have the following lemma for $x\in J(\theta)$.
\begin{lem}
\label{lem:BoundMixtureByStepFunc}
Under the assumption~\ref{assumption:1}, for each  $x \in J(\theta)$  
\[
  \sum_{m=1}^{M}
  \wtparam_{m} f_{m}(x; \locparam_{m}, \scparam_{m})
  \leq
  \sum_{t=1}^{T(\theta)}
  H(J_{t}(\theta)) \cdot \indicator_{J_{t}(\theta)}(x).
\]
\end{lem}

A density can be high only in a small region and we want to have some
explicit bound on the length $W(J_{t}(\theta))$ of $J_{t}(\theta)$
in terms of its height $H(J_{t}(\theta))$.
Let
\begin{equation}
   v_{2}
   \equiv 
   2\left(
     \frac{v_{1}}{\kappa_{0}}
    \right)^{\frac{1}{\beta}}
   \left(
    {v_{0}\cdot {(M + 1)}}
   \right)^{\tilde{\beta}}
  \quad , \quad 
  \xi(y) \equiv v_{2}\cdot \left(\frac{1}{y}\right)^{\tilde{\beta}}, \ y>0
  \ , 
  \label{eq:def:v_2_and_xi}
\end{equation}
where $v_0, v_1$ and $\beta$ are given in assumption~\ref{assumption:1}
and $\tilde \beta$ is defined in  (\ref{eq:def:beta_and_nu_and_v}).
\begin{lem}
\label{lem:BoundWidthOfInterval}
Under the assumption~\ref{assumption:1}, the length $W(J_{t}(\theta))$
of $J_{t}(\theta)$ for each $t$ is bounded as
 \begin{equation}
  W(J_{t}(\theta)) 
   \leq 
   v_{2} \cdot 
   \left(
    \frac{1}{{H(J_{t}(\theta))}}
   \right)^{\tilde{\beta}} 
   = \xi(H(J_{t}(\theta)))
   \ . 
   \nonumber 
 \end{equation}
\end{lem}
\Proof
In (\ref{eq:def:HJtn}) at any $x\in J_{t}(\theta)$, 
$H(J_{t}(\theta)) - \kappa_{0}$ 
consists of at most $M$ components.
Thus, for each $J_{t}(\theta)$, there exists at least one component $m=m_t$ such that 
\begin{eqnarray}
\frac{v_0}{\scparam_{m}}
  \geq
  \frac{1}{M}
  \cdot
  \left(
   H(J_{t}(\theta)) -{\kappa_{0}}
  \right) \; . 
  \nonumber 
\end{eqnarray}
Furthermore from (\ref{eq:kappa-c0-global-condition}) we have  
\begin{eqnarray}
  H(J_{t}(\theta)) - {\kappa_{0}}
  & > & 
  H(J_{t}(\theta)) - \frac{v_{0}}{c_{0}(M + 1)} 
  \nonumber \\
  & > & 
  H(J_{t}(\theta)) -  \frac{H(J_{t}(\theta))}{M + 1}
  = 
  \frac{M \cdot H(J_{t}(\theta))}{M + 1}
  \nonumber 
  .
\end{eqnarray}
Therefore we have
\[
  \nu(\scparam_{m}) 
=
  \left(
   \frac{v_{1}}{\kappa_{0}}
  \right)^{\frac{1}{\beta}}
  \scparam_{m}^{\tilde\beta}
\leq 
  \left(
   \frac{v_{1}}{\kappa_{0}}
  \right)^{\frac{1}{\beta}}
  \left(
   \frac{v_{0}{(M + 1)}}{{H(J_{t}(\theta))}}
  \right)
  ^{ \tilde\beta
  }  \ .
\]
This with $W(J_{t}(\theta)) \leq 2\nu(\scparam_{m})$ 
proves the lemma.
\qed
\bigskip 

So far we have been concerned with bounding the density at its peaks.
Now we consider bounding the tail of the true density $f(x;\theta_0)$.
Write 
$\bar{\locparam}_{0} \equiv \max(|\locparam_{01}|, \ldots, |\locparam_{0M}|)$
and
$
\theta_0=(\wtparam_{01},\locparam_{01},\scparam_{01}, \ldots, \wtparam_{0M},
\locparam_{0M}, \scparam_{0M}).
$
Let 
\begin{equation}
 u_0 \equiv \sup_x f(x;\theta_0) 
  \quad , \quad
  u_1 \equiv \max( u_0 \cdot (2 \bar{\locparam}_{0})^{\beta} , \; 
  2^{\beta} v_1 \sum_{m=1}^M \wtparam_{0m}
  \scparam_{0m}^{\beta-1}) 
  \; .
  \label{eq:def:u_0_and_u_1}
\end{equation}
\begin{lem}
\label{lem:true-tail}
Under the assumption~\ref{assumption:1}, 
the following inequality holds. 
\begin{eqnarray}
 f(x;\theta_{0})
  \leq
  \min{\{u_{0}, \; u_{1}\cdot \abs{x}^{-\beta}\}}, \quad \forall x\in \real. 
  \nonumber 
\end{eqnarray}
\end{lem}
\Proof
From assumption~\ref{assumption:1}
$$
u_0 \le \sum_{m=1}^M \wtparam_{0m} \frac{v_0}{\scparam_{0m}}.
$$
Then for $|x| \ge 2\bar{\locparam}_{0}$
$$
|x-\locparam_{0m}|^{-\beta} \le (|x| - \bar{\locparam}_{0})^{-\beta}
\le 2^\beta |x|^{-\beta} \ ,\quad  (m=1,\ldots,M).
$$
Therefore for $|x| \ge 2\bar{\locparam}_{0}$
$$
f(x;\theta_0) \le |x|^{-\beta} 2^\beta v_1 \sum_{m=1}^M \wtparam_{0m}
\scparam_{0m}^{\beta-1}
$$
and
\begin{eqnarray}
 f(x;\theta_{0})
  \leq
  \min{\{u_{0}, u_{1}\cdot \abs{x}^{-\beta}\}}, \quad \forall x\in \real.
  \nonumber 
\end{eqnarray}
\qed
\bigskip

Based on lemma \ref{lem:true-tail} we can bound
the behavior of the minimum and the maximum of the sample.
Let  $x_{1}, \ldots, x_{n}$ denote a random sample 
of size $n$ from $f(x;\theta_{0})$ and let 
\begin{eqnarray}
 x_{n,1}  \equiv  \min{\{x_{1}, \ldots, x_{n}\}}
  \quad , \quad 
 x_{n,n}  \equiv  \max{\{x_{1}, \ldots, x_{n}\}}.
  \nonumber 
\end{eqnarray}
The following lemma follows from the Borel-Cantelli lemma.
\begin{lem}
 \label{expanding-A:lem}
For any real constant $A_{0} > 0$ and $\zeta > 0$, 
define 
\begin{eqnarray}
 A_{n} \equiv A_{0} \cdot n^{\frac{2 + \zeta}{\beta - 1}}.
  \nonumber 
\end{eqnarray}
Then 
\begin{eqnarray}
 \prob
  \left(
   x_{n,1} < -A_{n} 
   \quad \rmor \quad 
   x_{n,n} > A_{n}
   \quad i.o.
  \right)
  = 0 \ .
  \nonumber 
\end{eqnarray}
\end{lem}
\Proof
By the Borel-Cantelli lemma and the Bonferroni inequality
it suffices to show that  
\[
 \sum_{n=1}^{\infty}
  \prob
  \left(
   x_{n,1} < -A_{n} 
  \right)
 <  
 \infty ,  \qquad
 \sum_{n=1}^{\infty}
 \prob
  \left(
   x_{n,n} > A_{n}
  \right)
  <  
 \infty .
\]
We consider the left tail. Let $F_{0}(x)$ denotes the distribution 
function of $f(x;\theta_{0})$.
Then 
\begin{eqnarray}
 \prob
  \left(
   x_{n,1} < -A_{n} 
  \right)
 & = & 
 1 - (1 - F_{0}(-A_{n}))^{n} ,
   \nonumber 
\end{eqnarray}
and 
\begin{eqnarray}
 F_{0}(-A_{n})
  & \leq & 
  \int_{-\infty}^{-A_{n}}
  u_{1} \cdot \abs{x}^{-\beta}
  \rmd x 
  = 
  \frac{u_{1}}{\beta - 1}
  \cdot A_{n}^{-\beta + 1} . 
   \nonumber 
\end{eqnarray}
By replacing $n$ by $n-n_0$ with a  sufficiently large $n_0$ if necessary,
we can assume without loss
of generality that 
\begin{eqnarray}
 \frac{u_{1}A_{0}^{-\beta + 1}}{\beta - 1} 
  \left(
   n^{\frac{2 + \zeta}{\beta - 1}}
  \right)^{-\beta + 1} < 1, \quad \forall n.
  \nonumber 
\end{eqnarray}
Then 
\begin{eqnarray}
 \log  (1 - F_{0}(-A_{n}))^{n}  
  & \geq & 
  \log{
   \left(
    1 - \frac{u_{1}}{\beta - 1}\cdot A_{n}^{-\beta + 1}
   \right)^{n}
  }
  \nonumber \\
  & = & 
   \log{
    \left(
     1 - \frac{u_{1}A_{0}^{-\beta + 1}}{\beta - 1} 
     \left(
      n^{\frac{2 + \zeta}{\beta - 1}}
     \right)^{-\beta + 1}
    \right)^{n}
   } .
   \nonumber 
\end{eqnarray}
Let 
$u_{2} \equiv u_{1}A_{0}^{-\beta + 1}/(\beta - 1)$
and we have 
\begin{eqnarray}
 \abslr{
  \log{
  (1 - F_{0}(-A_{n}))^{n}
  }}
 & \leq & 
  \abslr{
  \log{
  \left(1 - \frac{u_{2}}{n^{2 + \zeta}}\right)^{n}
  }} 
  \nonumber \\
 & = & 
 \frac{u_{2}}{n^{1 + \zeta}}
 \abslr{
 \log{
 \left(
  1 - \frac{u_{2}}{n^{2 + \zeta}}
 \right)^{\frac{n^{2 + \zeta}}{u_{2}}}
 }}
  \nonumber \\
 & = & 
  O(n^{-(1 + \zeta)}).
  \nonumber 
\end{eqnarray}
Hence there exists a  sufficiently large 
$N$ and $u_{3} > 0$
such that 
\begin{eqnarray}
  \abslr{
   \log{
  (1 - F_{0}(-A_{n}))^{n}
  }}
  \leq
  \frac{u_{3}}{n^{1 + \frac{\zeta}{2}}}
  \nonumber 
\end{eqnarray}
for all $n > N$.
This and 
$(1 - F_{0}(-A_{n}))^{n} \leq 1$ imply that  for $n > N$
\begin{eqnarray}
  \log{
  (1 - F_{0}(-A_{n}))^{n}
  }
  \geq
  -\frac{u_{3}}{n^{1 + \frac{\zeta}{2}}}.
  \nonumber 
\end{eqnarray}
Hence by $1-e^{-y} \leq y$, we have  for $n>N$
\begin{eqnarray}
 \prob(x_{n,1} < -A_{n})
  = 
 1 - (1 - F_{0}(-A_{n}))^{n}
  \leq 
 1 - \exp{\left(-\frac{u_{3}}{n^{1 + \frac{\zeta}{2}}}\right)}
 \leq
 \frac{u_{3}}{n^{1 + \frac{\zeta}{2}}}.
 \nonumber 
\end{eqnarray}
Therefore we obtain 
\begin{eqnarray}
 \sum_{n > N}
  \prob(x_{n,1} < -A_{n})
  =  
  \sum_{n > N} 1 - (1 - F_{0}(-A_{n}))^{n}
 \leq 
 \sum_{n} \frac{u_{3}}{n^{1 + \frac{\zeta}{2}}} 
 < \infty .
  \nonumber 
\end{eqnarray}
The case of  the right tail $\prob(x_{n,n} > A_{n})$ is also proved 
by the same argument.
\qed
\bigskip

Finally we consider subprobability densities in $\scrg_{\scrk}$.
For any positive real number $\rho$, let 
\begin{eqnarray}
 f_{\scrk}(x;\theta_{\scrk},\rho) 
 & \equiv & 
 \sup_{\dist(\theta_{\scrk}',\theta_{\scrk}) \leq \rho}f_{\scrk}(x;\theta_{\scrk}')
 \quad , \quad 
 (\theta_{\scrk}' \in  \bar{\Theta}_{\scrk}) .
 \label{eq:def:f_scrk_rho}
\end{eqnarray}
The following lemma follows from the bounded convergence theorem.
\begin{lem}
 \label{expectation-equality:lem}
 Let 
 ${\Gamma}_{\scrk}$ 
 denote any compact subset of $ \bar{\Theta}_{\scrk}$.
 For any real constant $\kappa \geq 0$ and 
 any point $\theta_{\scrk} \in {\Gamma}_{\scrk}$, 
 the following equality holds 
 under the assumption~\ref{assumption:1} and \ref{assumption:3}.
 \begin{eqnarray}
  \lim_{\rho \rightarrow 0}
   E_{0}[\log\{f_{\scrk}(x;\theta_{\scrk}, \rho) + \kappa\}]
    = 
    E_{0}[\log\{f_{\scrk}(x;\theta_{\scrk}) + \kappa\}] \ .
    \nonumber 
 \end{eqnarray}
\end{lem}
\Proof
We treat the case of $\kappa > 0$. 
The case of $\kappa = 0$ is almost the same as 
the proof of Lemma 2 in \cite{W1949}.
{}From assumption  \ref{assumption:3} we have 
\begin{eqnarray}
 \lim_{\rho \rightarrow 0}
  \log\{f_{\scrk}(x;\theta_{\scrk}, \rho) + \kappa\}
  = 
  \log\{f_{\scrk}(x;\theta_{\scrk}) + \kappa\}
  \qquad a.e.
  \nonumber 
\end{eqnarray}
Now 
${\Gamma}_{\scrk}$
is compact and $\kappa > 0$.
Hence by assumption  \ref{assumption:1}, 
$\log\{f_{\scrk}(x;\theta_{\scrk}, \rho) + \kappa\}$
is bounded.
Therefore 
\begin{eqnarray}
 \lim_{\rho \rightarrow 0}
  E_{0}[\log\{f_{\scrk}(x;\theta_{\scrk}, \rho) + \kappa\}]
  = 
  E_{0}[\log\{f_{\scrk}(x;\theta_{\scrk}) + \kappa\}]
  \nonumber 
\end{eqnarray}
by the bounded convergence theorem.
\qed
\bigskip

\subsection{Proof of theorem~\ref{kappa-lambda:thm}}
\label{sec:proof-kappa-lambda-thm}

In this section we prove theorem~\ref{kappa-lambda:thm} by
contradiction.  Fix arbitrary proper subset $\scrl$ of $\{1,\dots , M\}$.
It suffices to prove that (\ref{kappa-lambda:eq}) holds for 
all $g \in \scrg_{\scrl}$. 
Suppose that~(\ref{kappa-lambda:eq}) does not hold for some 
$\scrg_{\scrl}$. 
Then  for any $\lambda, \kappa > 0$, 
there exists $g \in \scrg_{\scrl}$ such that 
\begin{eqnarray}
 E_{0}[\log{\{g(x) + \kappa\}}] + \lambda \geq E_{0}[\log{f(x;\theta_{0})}].
  \nonumber 
\end{eqnarray}
Here, let $\{\lambda_{j}\}, \{\kappa_{j}\}$
be positive sequences which decrease 
to zero.
Then for each $\lambda_{j}, \kappa_{j} > 0$, 
there exists $g_{j} \in \scrg_{\scrl}$ such that 
\begin{eqnarray}
 E_{0}[\log{\{g_{j}(x) + \kappa_{j}\}}] + \lambda_{j}
  \geq E_{0}[\log{f(x;\theta_{0})}] \; .
  \nonumber 
\end{eqnarray}
It follows that
\begin{eqnarray}
 \liminf_{\seqnum \rightarrow \infty}
 E_{0}[\log{\{g_{j}(x) + \kappa_{j}\}}] + \lambda_{j}
  \geq E_{0}[\log{f(x;\theta_{0})}] \; .
  \label{eq:contradictory_assumption}
\end{eqnarray}
Now $g_{j}$ can be written as 
\begin{eqnarray}
 g_{j}(x) = f_{\scrl}(x; \theta_{\scrl}^{(\seqnum)}) \; .
  \nonumber 
\end{eqnarray}
Then the following lemma holds by compactification argument. 
\begin{lem}
\label{lem:proof-theor-refk}
 There exists a subsequence of 
 $\{\theta_{\scrl}^{(\seqnum)}\}_{\seqnum=1}^{\infty} \equiv \{\{\wtparam_{m}^{(\seqnum)}, \locparam_{m}^{(\seqnum)}, \scparam_{m}^{(\seqnum)} \mid m \in \scrl\}\}_{\seqnum=1}^{\infty}$ 
 and disjoint subsets 
 $
 \scrk_{{\scparam}\downarrow 0}, 
 \scrk_{{\scparam}\uparrow \infty}, 
 \scrk_{\abs{\locparam}\uparrow \infty}
 \subset \scrl 
 $ 
such that along the subsequence: 
\begin{align*}
&\scparam_{m}^{(\seqnum)} \rightarrow 0 \text{\quad for\quad} m \in
\scrk_{{\scparam}\downarrow 0}, \\
&\scparam_{m}^{(\seqnum)} \rightarrow \infty 
\text{\quad for\quad} m \in \scrk_{{\scparam}\uparrow \infty},\\
&\scparam_{m}^{(\seqnum)} \text{ converges to a finite value and }
\abs{\locparam_{m}^{(\seqnum)}} \rightarrow \infty \text{\ for\ }
m \in \scrk_{\abs{\locparam}\uparrow \infty},\\
&(\wtparam_{{m}}^{(\seqnum)}, \locparam_{{m}}^{(\seqnum)} ,
\scparam_{{m}}^{(\seqnum)}) \text{ converges to  a finite point } 
(\wtparam_{{m}}^{(\infty)}, \locparam_{{m}}^{(\infty)} , \scparam_{{m}}^{(\infty)}) 
 \text{ for } m \in \scrkrest,
\end{align*}
 where $\scrkrest \equiv \scrl 
 \backslash \{\scrk_{{\scparam}\downarrow 0} \cup \scrk_{{\scparam}\uparrow \infty} 
 \cup \scrk_{\abs{\locparam}\uparrow \infty}\}$.
\end{lem}
\Proof
Let 
\begin{eqnarray}
 {\locparam'}_{m}^{(\seqnum)} \equiv \arctan{(\locparam_{m}^{(\seqnum)})}, 
\qquad
 {\scparam'}_{m}^{(\seqnum)} \equiv \arctan{(\scparam_{m}^{(\seqnum)})}.
  \nonumber 
\end{eqnarray}
Then 
$\{\{\wtparam_{m}^{(\seqnum)}, {\locparam'}_{m}^{(\seqnum)}, {\scparam'}_{m}^{(\seqnum)} \mid m \in \scrl\}\}_{\seqnum=1}^{\infty}$
is regarded as a sequence in the following compact set. 
\begin{equation}
 0 \leq \wtparam_{m} \leq 1 
  \ , \ 
  \sum_{m = 1}^{L} \wtparam_{m} \leq 1 \ ,\ 
 -\frac{\pi}{2} \leq {\locparam'}_{m}^{(j)} \leq \frac{\pi}{2} \ , \ 
 0 \leq {\scparam'}_{m}^{(j)} \leq \frac{\pi}{2}.
 \label{extended-parameter-space:eq}
\end{equation}
Therefore there exists a subsequence of 
$\{\{\wtparam_{m}^{(j)}, {\locparam'}_{m}^{(j)}, {\scparam'}_{m}^{(j)} \mid m \in \scrl\}\}_{j=1}^{\infty}$ 
that converges to a point in the set (\ref{extended-parameter-space:eq}). 
Now we sort the elements in $\scrl$ according to their behaviors in this subsequence. 
First, we choose components such that 
$\scparam_{m}^{(j)} \rightarrow 0$, 
and add these $m$ to $\scrk_{{\scparam}\downarrow 0}$. 
Second, from the remainder, 
we choose components such that 
$\scparam_{m}^{(j)} \rightarrow \infty$, 
and add these $m$ to $\scrk_{{\scparam}\uparrow \infty}$.
Third, from the remainder, 
we choose components such that 
$\abs{\locparam_{m}^{(j)}} \rightarrow \infty$ 
and add these $m$ to $\scrk_{\abs{\locparam}\uparrow \infty}$. 
Finally, we choose 
the remaining components 
as $\scrkrest$.
\qed

\bigskip
From lemma~\ref{lem:proof-theor-refk}, we define $g_{\infty}$ as follows.  
\begin{eqnarray}
 g_{\infty}(x)
  \equiv 
  \sum_{m \in \scrkrest}
  {\wtparam}_{m}^{(\infty)} f_{m}(x; {\locparam}_{m}^{(\infty)}, {\scparam}_{m}^{(\infty)})
  \in \scrg_{\scrkrest} \; . 
  \nonumber 
\end{eqnarray}
For notational simplicity and without loss of generality, 
we replace the original sequence with this subsequence, because 
(\ref{eq:contradictory_assumption}) holds for this subsequence as well. 
Furthermore, by considering the sequence 
$\{\theta_{\scrl}^{(\seqnum)}\}_{\seqnum=\seqnum_{0}}^{\infty}$
where $j_{0}$ is sufficiently large and 
replacing $j$ by $j-j_0$ if necessary, 
we can assume without loss
of generality that there exist
sufficiently small real constants 
$\kappa_{0} > 0$ and $c_{0} > 0$ 
such that 
\begin{eqnarray}
 & &
  E_{0}[\log{f(x;\theta_{0})}]
  -
  E_{0}
  \left[
   \log{
   \left\{
    g_{\infty}(x) + 3\kappa_{0}
    \right\}
   }
  \right]
  > 
  0
  \;\; , \;\;
  {\kappa_{0}}
  < 
  \frac{v_0}{c_0(M + 1)} 
  \; , \nonumber \\
 & &
  \scparam_{m}^{(j)} < c_{0}
  \;\; (m \in \scrk_{{\scparam}\downarrow 0})
  \;\; , \;\; 
  \scparam_{m}^{(j)} > \frac{v_{0}}{\kappa_{0}}
  \;\  (m \in \scrk_{{\scparam}\uparrow \infty})
  \; , \nonumber \\
 & &
  c_{0} \leq \scparam_{m}^{(j)} \leq \frac{v_{0}}{\kappa_{0}} 
  \;\ (m \in \scrk_{\abs{\locparam}\uparrow \infty})
  \quad \text{for all }j \; . 
  \label{kappa-zero-and-c-zero:eq}
\end{eqnarray}
From lemma~\ref{lem:BoundComponentByStepFunc} and \ref{lem:BoundMixtureByStepFunc}, we have 
\begin{eqnarray}
 \lefteqn{
  E_{0}
  \left[
   \log{
   \left\{
    f_{\scrl}(x; {\theta}_{\scrl}^{(\seqnum)})
    + \kappa_{\seqnum}
   \right\}
   }
   + \lambda_{\seqnum}
  \right]
  }
  \nonumber \\
 & \leq & 
  \int 
  \indicator_{J({\theta}_{\scrl}^{(\seqnum)})}(x)
  \cdot
  \log{
  \left\{
   \sum_{t=1}^{T({\theta}_{\scrl}^{(\seqnum)})}
   H(J_{t}({\theta}_{\scrl}^{(\seqnum)}))\cdot \indicator_{J_{t}({\theta}_{\scrl}^{(\seqnum)})}(x)
   + \kappa_{\seqnum}
  \right\}
  }
  f(x;\theta_{0})
  \rmd x
  \nonumber \\
 & & + 
  \int 
  \indicator_{\real \backslash J({\theta}_{\scrl}^{(\seqnum)})}(x)
  \cdot
  \log{
  \left\{
   f_{\scrk_{{\scparam} \ndownarrow 0}}(x; \theta_{\scrk_{{\scparam} \ndownarrow 0}}^{(\seqnum)})
   +
   {\kappa_{0}}
   + \kappa_{\seqnum}
  \right\}
  }
  f(x;\theta_{0})
  \rmd x
   + \lambda_{\seqnum} \; , 
   \nonumber \\
  \label{inequality-of-expectation-of-g-infty:eq}
\end{eqnarray}
where 
$\scrk_{{\scparam} \ndownarrow 0} \equiv \scrl 
 \backslash \scrk_{{\scparam}\downarrow 0}$. 

Now we evaluate the first term 
on the right-hand side of (\ref{inequality-of-expectation-of-g-infty:eq}).
{}From lemma \ref{lem:BoundWidthOfInterval} 
\begin{eqnarray}
 \lefteqn{
  \int 
  \indicator_{J({\theta}_{\scrl}^{(\seqnum)})}(x)
  \cdot
  \log{
  \left\{
   \sum_{t=1}^{T({\theta}_{\scrl}^{(\seqnum)})}
   H(J_{t}({\theta}_{\scrl}^{(\seqnum)})) \cdot \indicator_{J_{t}({\theta}_{\scrl}^{(\seqnum)})}(x)
   + \kappa_{\seqnum}
  \right\}
  }
  f(x;\theta_{0})
  \rmd x
  } & & 
  \nonumber \\
 & \leq & 
  \sum_{t=1}^{T({\theta}_{\scrl}^{(\seqnum)})}
  W(J_{t}({\theta}_{\scrl}^{(\seqnum)}))
  \cdot
  \log{
  \left\{
   H(J_{t}({\theta}_{\scrl}^{(\seqnum)}))
   + \kappa_{\seqnum}
  \right\}
  }
  \cdot
  u_{0}
 \longrightarrow 0 
 \quad , \quad 
  (n \rightarrow \infty)\ , 
  \label{jack:eq}
\end{eqnarray}
where $u_{0} = \sup_x f(x;\theta_0)$ defined in (\ref{eq:def:u_0_and_u_1}). 
Next we evaluate the second term 
on the right-hand side of (\ref{inequality-of-expectation-of-g-infty:eq}).
Let 
\begin{eqnarray}
 A^{({\seqnum})}
  \equiv 
 \min_{m \in \scrk_{\abs{\locparam}\uparrow \infty}}
  \left\{
   \min{\{\abs{{\locparam}_{m}^{(\seqnum)} + \nu({\scparam}_{m}^{(\seqnum)})} , \abs{{\locparam}_{m}^{(\seqnum)} - \nu({\scparam}_{m}^{(\seqnum)})}\}}
  \right\}. 
  \label{eq:def:A_hat_seqnum}
\end{eqnarray}
Then
\begin{eqnarray}
& & 
f_{\scrk_{{\scparam} \ndownarrow 0}}(x; \theta_{\scrk_{{\scparam} \ndownarrow 0}}^{(\seqnum)})
 = 
f_{\scrk_{{\scparam} \uparrow \infty}}
(x;\theta_{\scrk_{{\scparam} \uparrow \infty}}^{(\seqnum)})
+ 
f_{\scrk_{\abs{\locparam} \uparrow \infty}}
(x;\theta_{\scrk_{\abs{\locparam} \uparrow \infty}}^{(\seqnum)})
+ 
f_{\scrkrest}(x; \theta_{\scrkrest}^{(\seqnum)})
\; , 
\nonumber \\ 
& & f_{\scrk_{{\scparam} \uparrow \infty}}(x;\theta_{\scrk_{{\scparam} \uparrow \infty}}^{(\seqnum)}) \leq \kappa_{0}
\quad \text{ for all } x 
 \; , \nonumber \\ 
& & f_{\scrk_{{\scparam} \uparrow \infty}}
(x;\theta_{\scrk_{{\scparam} \uparrow \infty}}^{(\seqnum)})
 + f_{\scrk_{\abs{\locparam} \uparrow \infty}}
(x;\theta_{\scrk_{\abs{\locparam} \uparrow \infty}}^{(\seqnum)})
 \leq \kappa_{0} 
\quad \text{ for } x \in [-A^{({\seqnum})},A^{({\seqnum})}] \backslash J({\theta}_{\scrl}^{(\seqnum)}) \; . 
\nonumber 
\end{eqnarray}
Therefore the following inequality holds.
\begin{eqnarray}
 \lefteqn{
  \int 
  \indicator_{\real \backslash J({\theta}_{\scrl}^{(\seqnum)})}(x)
  \cdot
  \log{
  \left\{
   f_{\scrk_{{\scparam} \ndownarrow 0}}(x; \theta_{\scrk_{{\scparam} \ndownarrow 0}}^{(\seqnum)})
   +
   {\kappa_{0}}
   + \kappa_{\seqnum}
  \right\}
  }
  f(x;\theta_{0})
  \rmd x
  } & & 
  \nonumber \\
 & \leq & 
  \int 
  \indicator_{[-A^{({\seqnum})},A^{({\seqnum})}] \backslash J({\theta}_{\scrl}^{(\seqnum)})}(x)
  \cdot
  \log{
  \left\{
   f_{\scrkrest}(x; \theta_{\scrkrest}^{(\seqnum)})
   +
   {2\kappa_{0}}
   + \kappa_{\seqnum}
  \right\}
  }
  f(x;\theta_{0})
  \rmd x 
  \nonumber \\
 & & + 
  \int 
  \indicator_{
  \{(-\infty, -A^{({\seqnum})}) \bigcup (A^{({\seqnum})}, \infty)\}
  \backslash J({\theta}_{\scrl}^{(\seqnum)})}(x)
  \cdot
  \log{
  \left\{
   f_{\scrk_{\abs{\locparam}\uparrow \infty}}(x; \theta_{\scrk_{\abs{\locparam}\uparrow \infty}}^{(\seqnum)})
  \right.
  }
  \nonumber \\
 & & \qquad\qquad
  { + 
  \left.
   f_{\scrkrest}(x; \theta_{\scrkrest}^{(\seqnum)})
   +
   2{\kappa_{0}}
   + \kappa_{\seqnum}
  \right\}
  }
  f(x;\theta_{0})
  \rmd x 
  \nonumber   \\
 &\equiv& I_1^{(j)} + I_2^{(j)} \qquad{\rm(say)} .
  \label{tree-of-beans:eq}
\end{eqnarray}
By the bounded convergence theorem, 
we obtain 
\begin{eqnarray}
I_1^{(j)}
  \rightarrow 
  \int 
  \log{
  \left\{
   g_{\infty}(x)
   +
   {2\kappa_{0}}
  \right\}
  }
  f(x;\theta_{0})
  \rmd x
  \quad , \quad 
  I_2^{(j)}
   \rightarrow 
  0. 
  \label{fatboy-slim-2:eq}
\end{eqnarray} 

{}From (\ref{inequality-of-expectation-of-g-infty:eq}),
(\ref{jack:eq}), (\ref{tree-of-beans:eq}), 
(\ref{fatboy-slim-2:eq}), 
we have 
\begin{eqnarray}
 E_{0}[\log f(x;\theta_{0})]
  \leq
 \limsup_{j \rightarrow \infty}
  E_{0}[\log\{g_{j}(x) + \kappa_{j}\}] + \lambda_{j}
  \leq 
  E_{0}
  \left[
   \log{
   \left\{
    g_{\infty}(x)
    +
    {2\kappa_{0}}
   \right\}
   }
  \right].
  \nonumber 
\end{eqnarray}
This is a contradiction to (\ref{kappa-zero-and-c-zero:eq}).
This completes the proof of theorem \ref{kappa-lambda:thm}.

\subsection{Proof of the main theorem}
\label{sec:proof-main-thm}

We choose real constants $\kappa$ and $\lambda$ to satisfy (\ref{kappa-lambda:eq}) 
by using theorem \ref{kappa-lambda:thm}. 
Having chosen these constants, 
from now on we follow the line of the proof 
in \cite{TT2003-20}, although the details of the proof here is much more
complicated. 
For the sake of readability we divide our proof into further sections.

\subsubsection{Setting up constants}
For $\kappa, \lambda$
satisfying~(\ref{kappa-lambda:eq}), let $\kappa_{0}, \lambda_{0}$ be
real constants such that 
\begin{equation}
 0 < 4\kappa_{0} \leq \kappa \quad , \quad 0 < 4\lambda_{0} \leq \lambda \; .
  \nonumber 
\end{equation}
Note that $4\kappa_{0}, 4\lambda_{0}$ also satisfy~(\ref{kappa-lambda:eq}).
Define 
\begin{eqnarray}
 B \equiv \frac{v_{0}}{\kappa_{0}} 
  > 
  \max{\{\scparam_{01}, \ldots, \scparam_{0M}\}}.
  \label{eq:def:B}
\end{eqnarray}
If $\scparam_{m} \geq B$, then the density of the $m$-th component is almost flat 
and makes little contribution to the likelihood. 
In section \ref{sec:part-param-space}, we partition the parameter space 
according to this property. 

Because $\{c_n\}$ is decreasing
to zero, by replacing $c_0$ by some $c_n$ if necessary, 
we can assume without loss
of generality that $c_0$ is sufficiently small to satisfy the
following conditions, 
\begin{eqnarray}
 & & (v_0/c_0)^{\tilde{\beta}}
  > e , 
  \nonumber \\
 & & c_{0}  <  \min{\{\scparam_{01}, \ldots, \scparam_{0M}\}} , 
  \nonumber \\
 & & 3M \cdot u_{0} \cdot 2\nu(c_{0}) \cdot \abslr{\log{ {\kappa_{0}} }}
   < 
  {\lambda_{0}} ,
   \nonumber 
  \\ 
 & &
 3 \cdot 2M \cdot u_{0} \cdot
  \xi(v_0/c_{0})
  \cdot \log(v_0/c_0)
   < 
  {\lambda_{0}}
  \; , 
   \nonumber 
  \\
  & & {\kappa_{0}}
  < 
  \frac{v_0}{c_0(M + 1)}
   \quad , 
   \label{eq:condition:proof:mainthm}
\end{eqnarray}
where $\tilde{\beta}$, $\nu(\cdot)$
and $\xi(\cdot)$ are defined in 
(\ref{eq:def:beta_and_nu_and_v}) and (\ref{eq:def:v_2_and_xi}). 

For any subset $V\subset \real$, let $P_0(V)$ denote the probability of
$V$ under the true density
\begin{eqnarray}
 P_0(V) \equiv \int_{V} f(x;\theta_{0}) \rmd x \; . 
  \label{eq:def:P_0}
\end{eqnarray}
Let $A_{0} > 0$ be a positive constant which satisfies  
\begin{eqnarray}
  P_0(\scra_{0})
  \cdot
  \log{
  \left(
   \frac{v_0/c_0 + 2\kappa_{0}}{3\kappa_{0}}
  \right)
  }
  < {\lambda_{0}} , 
  \label{A-zero:condition:eq}
\end{eqnarray}
where
\begin{equation}
 \scra_{0} \equiv (-\infty, -A_{0}] \cup [A_{0},\infty).
  \label{eq:def:scra_0}
\end{equation}
Let
$
 A_{n} \equiv A_{0} \cdot n^{\frac{2 + \zeta}{\beta - 1}} 
  \nonumber 
$
as in lemma \ref{expanding-A:lem}.
Define a subset $\Theta_{n}'$ of $\Theta_{n}$ in
theorem~\ref{main-thm:thm} by 
\begin{eqnarray}
  \Theta_{n}'
  \equiv
  \{
  \theta \in \Theta_{n}
  \mid
  \exists m \; {\rm s.t.} \; 
  c_{n} \leq \scparam_{m} \leq c_{0}
  \; \mathrm{or} \; 
  \abs{\locparam_{m}} > A_{0}
  \} \subset \Theta_{n} \ ,  
  \nonumber
\end{eqnarray}
and let
\begin{eqnarray}
 \Gamma_{0}
  \equiv
  \{
  \theta \in \Theta
  \mid
  c_{0} \leq \scparam_{m} \leq B
  \; , \;
  \abs{\locparam_{m}} \leq A_{0}
  ,\  (m=1,\ldots,M)\ 
  \}
  \; .
  \nonumber
\end{eqnarray}
Note that $\Theta_0 \subset \Gamma_0$, where $\Theta_0$ is the set of
true parameters.

\subsubsection{Partitioning the parameter space}
\label{sec:part-param-space}

In view of theorems in \cite{W1949}, \cite{R1981}, 
for the strong consistency of MLE on $\Theta_n$ 
under assumption~\ref{assumption:1}, \ref{assumption:2}, \ref{assumption:3} and \ref{assumption:4}, 
it suffices to prove that 
$$
 \lim_{n \rightarrow \infty}
  \frac{
   \sup_{\theta \in \Gamma \cup \Theta_{n}'}
   \prod_{i=1}^{n} f(x_{i};\theta)
   }
   { \prod_{i=1}^{n} f(x_{i};\theta_{0}) } = 0, \quad a.e. 
$$
for all closed $\Gamma\subset \Gamma_0$ not intersecting $\Theta_0$. 
Note that for all $\Gamma$ and $\{x_{i}\}_{i=1}^{n}$, 
\begin{eqnarray}
 {
  \sup_{\theta \in \Gamma \cup \Theta_{n}'}
  \prod_{i=1}^{n} f(x_{i};\theta)
  }
  = 
  \max
  \left\{
  {
  \sup_{\theta \in \Gamma}
  \prod_{i=1}^{n} f(x_{i};\theta)
  }
  \; , \; 
  {
  \sup_{\theta \in \Theta_{n}'}
  \prod_{i=1}^{n} f(x_{i};\theta)
  }
  \right\}
  \; . \;
  \nonumber
\end{eqnarray}
Furthermore
$$
   \lim_{n \rightarrow \infty}
  \frac{
   \sup_{\theta \in \Gamma}
   \prod_{i=1}^{n} f(x_{i};\theta)
   }
   { \prod_{i=1}^{n} f(x_{i};\theta_{0}) }
  = 0, \quad a.e.
$$
holds by 
theorems in \cite{W1949}, \cite{R1981}. 
Therefore it suffices to prove 
\begin{eqnarray}
   \lim_{n \rightarrow \infty}
   \frac{
   \sup_{\theta \in \Theta_{n}'}
   \prod_{i=1}^{n} f(x_{i};\theta)
   }
   { \prod_{i=1}^{n} f(x_{i};\theta_{0}) }
   = 0, \quad a.e.
   \nonumber 
\end{eqnarray}

Note that in the argument above the supremum of the likelihood
function over $\Gamma \cup \Theta_{n}'$ is considered separately for $\Gamma$
and $\Theta_n'$.  $\Gamma$ and $\Theta_n'$ form a covering of $\Gamma\cup
\Theta_n'$.  In our proof, we consider finer and finer finite coverings of
$\Theta_n'$.  As above, it suffices to prove that the ratio of the
supremum of the likelihood over each member of the covering to the
likelihood at $\theta_0$ converges to zero almost everywhere. 

Let $\theta \in \Theta_{n}'$.
Let 
$
 \scrk_{{\scparam} \leq c_{0}}, 
 \scrk_{{\scparam}\geq B}, 
 \scrk_{\abs{\locparam} \geq A_{0}} 
$
represent disjoint subsets of $\{1,\dots,M\}$ and 
define 
\begin{equation}
\scrkrest \equiv \{1, \dots, M\} 
\backslash 
\{\scrk_{{\scparam} \leq c_{0}} \cup \scrk_{{\scparam}\geq B} \cup \scrk_{\abs{\locparam} \geq A_{0}}\}
 . 
 \nonumber 
\end{equation}
For any given 
$ 
 \scrk_{{\scparam} \leq c_{0}}, 
 \scrk_{{\scparam}\geq B}, 
 \scrk_{\abs{\locparam} \geq A_{0}} 
$, 
we define  a subset of $\Theta_n'$ by
\begin{eqnarray}
 \Theta_{n, \scrk}'
  & \equiv &
  \{
  \theta \in \Theta_{n}'
   \mid 
  \scparam_{m} \leq c_{0},(m \in \scrk_{{\scparam} \leq c_{0}})
  \; ; \;
  \scparam_{m} \geq B,(m \in \scrk_{{\scparam}\geq B})
  \; ; \;
  \nonumber \\
 && 
  \qquad \qquad \quad 
  c_{0} < \scparam_{m} < B,\; \abs{\locparam_{m}}\geq A_{0},\; (m \in \scrk_{\abs{\locparam} \geq A_{0}})
  \; ; \;
  \nonumber \\
 && 
  \qquad \qquad \quad 
  c_{0} < \scparam_{m} < B,\; \abs{\locparam_{m}}<A_{0},\; (m \in \scrkrest)
  \}. 
  \label{eq:def:Theta_n_scrk}
\end{eqnarray}
As above, it suffices to prove that for each choice of disjoint subsets 
$
 \scrk_{{\scparam} \leq c_{0}}, 
 \scrk_{{\scparam}\geq B}, 
 \scrk_{\abs{\locparam} \geq A_{0}}, 
$ 
\begin{eqnarray}
  \lim_{n \rightarrow \infty}
  \frac{
  \sup_{\theta \in \Theta_{n,\scrk}'}
  \prod_{i=1}^{n} f(x_{i};\theta)
  }
  { \prod_{i=1}^{n} f(x_{i};\theta_{0}) }
  = 0, \quad a.e. 
  \nonumber 
\end{eqnarray}
We fix 
$
 \scrk_{{\scparam} \leq c_{0}}, 
 \scrk_{{\scparam}\geq B}, 
 \scrk_{\abs{\locparam} \geq A_{0}}, 
$ 
from now on.

Next we consider coverings of $\bar{\Theta}_{\scrkrest}$. 
Recall that $\bar{\Theta}_{\scrk}$, $f_{\scrk}(x;\theta_{\scrk})$ and $f_{\scrk}(x;\theta_{\scrk}, \rho)$ 
are defined in (\ref{eq:def:barTheta_scrk}), (\ref{eq:def:f_scrk}) and (\ref{eq:def:f_scrk_rho}). 
The following lemma follows from lemma \ref{expectation-equality:lem} and compactness of $\bar{\Theta}_{\scrkrest}$. 
\begin{lem}
\label{WaldThm1Type:lem}
Let $\scrb(\theta,\rho(\theta) )$ denote the open ball
with center $\theta$ and radius $\rho(\theta)$.
Then $\bar{\Theta}_{\scrkrest}$ can be covered by 
a finite number of balls 
$\scrb(\theta_{\scrkrest}^{(1)}, \rho(\theta_{\scrkrest}^{(1)}))
, \ldots, \scrb(\theta_{\scrkrest}^{(S)}, \rho(\theta_{\scrkrest}^{(S)})) 
$
such that 
\begin{eqnarray}
 E_{0}[\log{\{f_{\scrkrest}(x;\theta_{\scrkrest}^{(s)},\rho(\theta_{\scrkrest}^{(s)})) + \kappa_{0}\}}]
  + \lambda_{0}
  <
  E_{0}[\log{f(x;\theta_0)}]
  \; , \quad (s=1,\ldots,S) \; . 
  \nonumber 
\end{eqnarray}
\end{lem}
\Proof
{}From lemma \ref{expectation-equality:lem} we have
\begin{eqnarray}
 \lim_{\rho \rightarrow 0}
  E_{0}
  \left[
   \log{\{f_{\scrkrest}(x; \theta_{\scrkrest}, \rho) + \kappa_{0}\}}
  \right]
  = 
  E_{0}
  \left[
   \log{\{f_{\scrkrest}(x; \theta_{\scrkrest}) + \kappa_{0}\}}
  \right].
  \nonumber 
\end{eqnarray}
For each $\theta_{\scrkrest} \in \bar{\Theta}_{\scrkrest}$
\begin{eqnarray}
  E_{0}
  \left[
   \log{
   \left\{ 
    f_{\scrkrest}(x; \theta_{\scrkrest}) + \kappa_{0}
   \right\}
   }
  \right]
  +
  \lambda_{0}
  <
  E_{0}[\log{f(x;\theta_{0})}] 
  \nonumber 
\end{eqnarray}
holds.
Therefore for each 
$\theta_{\scrkrest} \in \bar{\Theta}_{\scrkrest}$, 
there exists a radius 
$\rho(\theta_{\scrkrest}) > 0$ such that 
\begin{eqnarray}
 E_{0}[\log{\{f_{\scrkrest}(x;\theta_{\scrkrest},\rho(\theta_{\scrkrest})) + \kappa_{0}\}}]
  + \lambda_{0}
  <
  E_{0}[\log{f(x;\theta_0)}].
  \nonumber 
\end{eqnarray}
Since 
\begin{eqnarray}
 \bar{\Theta}_{\scrkrest}
  \subset 
  \bigcup_{\theta_{\scrkrest} \in \bar{\Theta}_{\scrkrest}}
  \scrb(\theta_{\scrkrest},\rho(\theta_{\scrkrest}))
  \nonumber 
\end{eqnarray}
and the compactness of $\bar{\Theta}_{\scrkrest}$,
there exists a finite number of balls 
$\scrb(\theta_{\scrkrest}^{(1)}, \rho(\theta_{\scrkrest}^{(1)})),\ldots,$
\\
$\scrb(\theta_{\scrkrest}^{(S)}, \rho(\theta_{\scrkrest}^{(S)}))$
which cover $\bar{\Theta}_{\scrkrest}$.
\qed 

\bigskip
Based on lemma \ref{WaldThm1Type:lem} we partition $\Theta_{n,\scrk}'$.
Recall that we denote by $\theta_{\scrk}$ the subvector of $\theta \in \Theta$ 
consisting of the components in $\scrk$.
Define a subset of $\Theta_{n,\scrk}'$ by
\begin{eqnarray}
 \Theta_{n,\scrk ,s}'
  \equiv 
  \{
  \theta \in \Theta_{n,\scrk}'
  \mid 
  \theta_{\scrkrest} \in \scrb(\theta_{\scrkrest}^{(s)},\rho(\theta_{\scrkrest}^{(s)}))
  \}.
  \label{eq:def:Theta_scrk_n_s}
\end{eqnarray}
Then $\Theta_{n,\scrk}'$ is covered by
$\Theta_{n,\scrk,1}', \ldots, \Theta_{n,\scrk,S}'$ :
$$
\Theta_{n,\scrk}' = \bigcup_{s=1}^S \Theta_{n,\scrk,s}' \   .
$$
Again it suffices to prove
that for each choice of 
$
 \scrk_{{\scparam} \leq c_{0}}, 
 \scrk_{{\scparam}\geq B}, 
 \scrk_{\abs{\locparam} \geq A_{0}}, 
 s 
$ 
\begin{eqnarray}
   \lim_{n \rightarrow \infty}
   \frac{
   \sup_{\theta \in \Theta_{n,\scrk,s}'}
   \prod_{i=1}^{n} f(x_{i};\theta)
   }
   { \prod_{i=1}^{n} f(x_{i};\theta_{0}) }
   = 0, \quad a.e.
  \label{goal3:LikelihoodRatio:eq}
\end{eqnarray}
We fix 
$
 \scrk_{{\scparam} \leq c_{0}}, 
 \scrk_{{\scparam}\geq B}, 
 \scrk_{\abs{\locparam} \geq A_{0}}$ and 
$s$ 
from now on.
Because 
$$
\lim_{n\rightarrow\infty} \frac{1}{n} \sum_{i=1}^{n} 
\log{f(x_{i};\theta_0)} =
E_0[\log f(x;\theta_0)], \quad a.e.
$$
(\ref{goal3:LikelihoodRatio:eq}) is implied by
\begin{eqnarray}
\limsup_{n \rightarrow \infty} \frac{1}{n} \sup_{\theta \in \Theta_{n,\scrk,s}'}
   \sum_{i=1}^{n} \log{f(x_{i};\theta)}
  < 
   E_0[\log f(x;\theta_0)], 
  \quad a.e.
  \label{goal:eq}
\end{eqnarray}
Therefore it suffices to prove (\ref{goal:eq}), which is a new
intermediate goal of our proof hereafter.

\subsubsection{Bounding the likelihood by four terms}
In this section we bound the likelihood function by four terms
depending on the positions of the observations $x_1, \dots, x_n$.
Let $R_{n}(V)$ denote the number of observations
which belong to a set $V \subset \real$.
\begin{lem}
 \label{lem:fourterms}
For $\theta \in \Theta_{n,\scrk,s}'$
\begin{eqnarray}
 \frac{1}{n}\sum_{i=1}^{n}\log{f(x_{i};\theta)}
 & \leq & 
  \frac{1}{n}\sum_{i=1}^{n}
  \log{
  \left\{
   f_{\scrkrest}(x_{i};\theta_{\scrkrest},\rho(\theta_{\scrkrest})) + 3{\kappa_{0}}
  \right\}
  }
  \nonumber \\
 & & +
  \frac{1}{n}R_{n}(\scra_{0}) \cdot
  \log{
  \left(
   \frac{M v_0/c_0 + 2\kappa_{0}}{3\kappa_{0}}
  \right)
  } 
  \nonumber \\
 & & +
  \frac{1}{n}R_n(J(\theta))\cdot 
(
   -\log\kappa_{0} 
)
  + 
  \frac{1}{n}\sum_{x_{i} \in J(\theta)}
  \log{f(x_{i};\theta)}
  \; . 
  \label{loglikelihood:lem:eq}
\end{eqnarray}
\end{lem}
\Proof 
Let $\scrk_{{\scparam} > c_{0}} 
= \{1,\dots,M\}\backslash\scrk_{{\scparam} \leq c_{0}}$ 
and 
$\scrkbbound = \{1,\dots,M\}\backslash\{\scrk_{{\scparam} \leq c_{0}}\cup\scrk_{{\scparam} \geq B}\}$.
For $x \not\in J(\theta)$, 
$f(x;\theta)\leq f_{\scrk_{{\scparam} > c_{0}}}
(x;\theta_{\scrk_{{\scparam} > c_{0}}}) 
+ \kappa_{0}$ 
holds.
Therefore
\begin{eqnarray}
 \frac{1}{n}\sum_{i=1}^{n}\log{f(x_{i};\theta)}
  & \leq & \frac{1}{n} \sum_{x_i \in J(\theta)} \log{f(x_{i};\theta)}
  + \frac{1}{n} \sum_{x_i \not \in J(\theta)} 
  \log\left\{ 
       f_{\scrk_{{\scparam} > c_{0}}}(x;\theta_{\scrk_{{\scparam} > c_{0}}}) + {\kappa_{0}} 
      \right\}
  \nonumber \\
 &= &
    \frac{1}{n}\sum_{i=1}^{n}
  \log{
  \left\{
   f_{\scrk_{{\scparam} > c_{0}}}(x;\theta_{\scrk_{{\scparam} > c_{0}}}) + {\kappa_{0}} 
  \right\}
  }
  \nonumber \\
 && \qquad  +
  \frac{1}{n}\sum_{x_{i} \in J(\theta)}
  \left[
   \log
   f(x_i ; \theta)
   -
   \log
   \left\{
    f_{\scrk_{{\scparam} > c_{0}}}(x;\theta_{\scrk_{{\scparam} > c_{0}}}) + {\kappa_{0}} 
   \right\}
  \right]
\label{eq:lemma9proof1}
\end{eqnarray}
Consider the second term on the right-hand side. We have
\begin{align*}
&  \frac{1}{n}\sum_{x_{i} \in J(\theta)}
  \left[
   \log
   f(x_i ; \theta)
   -
   \log
   \left\{
    f_{\scrk_{{\scparam} > c_{0}}}(x;\theta_{\scrk_{{\scparam} > c_{0}}}) + {\kappa_{0}} 
   \right\}
  \right]\\
& \qquad\qquad \le 
   \frac{1}{n}\sum_{x_{i} \in J(\theta)}
    \log
    f(x_i ; \theta)
   -
   \frac{1}{n}R_n(J(\theta)) \cdot \log{\kappa_{0}}
   \; . 
\end{align*}
This takes care of the third and the fourth term of 
  (\ref{loglikelihood:lem:eq}).

Now consider the first term on the right-hand side of
(\ref{eq:lemma9proof1}).
Note that
\begin{eqnarray}
    \frac{1}{n}\sum_{i=1}^{n}
  \log{
  \left\{
   f_{\scrk_{{\scparam} > c_{0}}}(x;\theta_{\scrk_{{\scparam} > c_{0}}}) + {\kappa_{0}} 
  \right\}
  }
   &\leq &
  \frac{1}{n}\sum_{i=1}^{n}
  \log{
  \left\{
     f_{\scrkbbound}(x_{i};\theta_{\scrkbbound}) + 2{\kappa_{0}}
  \right\}
  }
  \nonumber 
\end{eqnarray}
For $x\notin  \scra_{0}$ 
$$
 f_{\scrk_{\abs{\locparam} \geq A_{0}}}(x;\theta_{\scrk_{\abs{\locparam} \geq A_{0}}})
  \leq
 {\kappa_{0}}.
$$
Therefore we obtain 
\begin{eqnarray}
 \lefteqn{
 \frac{1}{n}\sum_{i=1}^{n}
  \log{
  \left\{
     f_{\scrkbbound}(x_{i};\theta_{\scrkbbound})
   + 2{\kappa_{0}}
  \right\}
  } }
  \nonumber \\
  & = & 
  \frac{1}{n}\sum_{x_{i} \notin \scra_{0}}
  \log{
  \left\{
   f_{\scrkbbound}(x_{i};\theta_{\scrkbbound})
   + 2{\kappa_{0}}
  \right\}
  }
  + 
  \frac{1}{n}\sum_{x_{i} \in \scra_{0}}
  \log{
  \left\{
   f_{\scrkbbound}(x_{i};\theta_{\scrkbbound})
   + 2{\kappa_{0}}
  \right\}
  }
  \nonumber \\
 & \leq &
    \frac{1}{n}\sum_{x_{i} \notin \scra_{0}}
  \log{
  \left\{
   f_{\scrkrest}(x_{i};\theta_{\scrkrest})
   + 3{\kappa_{0}}
  \right\}
  }
  + 
  \frac{1}{n}\sum_{x_{i} \in \scra_{0}}
  \log{
  \left\{
   f_{\scrkbbound}(x_{i};\theta_{\scrkbbound})
   + 2{\kappa_{0}}
  \right\}
  }
  \nonumber \\
 & = & 
  \frac{1}{n}\sum_{i=1}^{n}
  \log{
  \left\{
   f_{\scrkrest}(x_{i};\theta_{\scrkrest})
   + 3{\kappa_{0}}
  \right\}
  }
  \nonumber \\
 & & \qquad \qquad 
  + 
  \frac{1}{n}\sum_{x_{i} \in \scra_{0}}
  \left[
   \log{
   \left\{
    f_{\scrkbbound}(x_{i};\theta_{\scrkbbound})
    + 2{\kappa_{0}}
   \right\}
   }
  - 
   \log{
   \left\{
   f_{\scrkrest}(x_{i};\theta_{\scrkrest})
    + 3{\kappa_{0}}
   \right\}
   }
  \right]
  \nonumber  \\ 
  \label{eq:totyu}
\end{eqnarray}
Note that $f_{\scrkbbound}(x;\theta_{\scrkbbound}) \leq v_{0}/c_{0}$
from lemma~\ref{lem:BoundComponentByStepFunc}. Therefore 
\begin{eqnarray}
 \lefteqn{
 \text{The r.h.s of (\ref{eq:totyu})} 
 } & & \nonumber \\
 &\leq &
  \frac{1}{n}\sum_{i=1}^{n}
  \log{
  \left\{
   f_{\scrkrest}(x_{i};\theta_{\scrkrest})
   + 3{\kappa_{0}}
  \right\}
  }
  + 
  \frac{1}{n}\sum_{x_{i} \in \scra_{0}}
  \left[
   \log{
   \left\{
    v_0/c_0
    + 2{\kappa_{0}}
   \right\}
   }
   - 
   \log{
   3\kappa_{0}
   }
  \right]
  \nonumber \\
 &\leq &
  \frac{1}{n}\sum_{i=1}^{n}
  \log{
  \left\{
     f_{\scrkrest}(x_{i};\theta_{\scrkrest},\rho(\theta_{\scrkrest})) + 3{\kappa_{0}}
  \right\}
  }
  + 
  \frac{1}{n} R_{n}(\scra_{0}) \cdot 
  \log{
  \left(
   \frac{v_0/c_0 + 2{\kappa_{0}}}{3{\kappa_{0}}}
  \right)
  }.
  \nonumber 
\end{eqnarray}
This takes care of the first and the second term of 
(\ref{loglikelihood:lem:eq}).
\qed

\bigskip 
From lemma~\ref{WaldThm1Type:lem} and the strong law of large
numbers the first term on the right hand side of (\ref{loglikelihood:lem:eq}) 
converges to the expectation of a density which has less than $M$ components 
and the expectation is less than that of the true density by theorem \ref{kappa-lambda:thm}. 
The second term 
converges to a small value because the relative frequency on $\scra_{0}$ 
is very small. 
The third term 
also converges to small value because the relative frequency on $J(\theta)$ 
is very small. 
The fourth term 
is somewhat complicated. 
The component in $\scrk_{\scparam \le c_{0}}$ may have high peaks. 
However the widths of the peaks are very narrow and the relative frequency on the interval is very small.
Hence the fourth term makes little contribution to the likelihood. 
Therefore the mean log likelihood (the left hand side of (\ref{loglikelihood:lem:eq})) converges to 
a value which is less than that of the true density. 
In the following we consider the details. 

The first term and the second term are easy.  

\bigskip
\noindent
{\bf The first term}: By lemma~\ref{WaldThm1Type:lem} and the strong law of large
numbers we have
\begin{eqnarray}
\label{bounding-barf:eq}
\lim_{n\rightarrow\infty}  \frac{1}{n}\sum_{i=1}^{n}
  \log{
  \left\{
     f_{\scrkrest}(x_{i};\theta_{\scrkrest},\rho(\theta_{\scrkrest})) + 4{\kappa_{0}}
  \right\}
  }
 < 
 E_0[\log f(x;\theta_0)] - 4\lambda_{0} , \quad a.e.
\end{eqnarray}

\bigskip
\noindent
{\bf The second term}:
By (\ref{A-zero:condition:eq}) and the strong law of large
numbers we have
\begin{eqnarray}
 \label{bounding-A-zero-term:eq}
  \lim_{n\rightarrow\infty} 
  \frac{1}{n}R_{n}(\scra_{0}) \cdot
  \log{
  \left(
   \frac{v_0/c_0 + 2\kappa_{0}}{3\kappa_{0}}
  \right)
  } 
  < 
  {\lambda_{0}}
  , \quad a.e.
\end{eqnarray}

Note that we have $-4\lambda_0$ from the first term and $\lambda_0$
from the second term. In the rest of our proof we show that
both the third term and the fourth term can be bounded by $\lambda_0$.

\subsubsection{Bounding the third term}

The third term can be bounded by dividing the interval $[-A_n, A_n]$
into short intervals of length $2\nu(c_0)$.
\begin{lem}
\label{boundedRJ:lem}
 \begin{eqnarray}
 \limsup_{n \rightarrow \infty}
  \sup_{\theta \in \Theta_{n,\scrk,s}'}
   \frac{1}{n} R_n(J(\theta))
   \leq 3M \cdot u_{0} \cdot 2\nu(c_{0}), 
  \quad a.e.
  \nonumber 
\end{eqnarray}
\end{lem}
\Proof 
\begin{figure}[htbp]
 \begin{center}
  \includegraphics{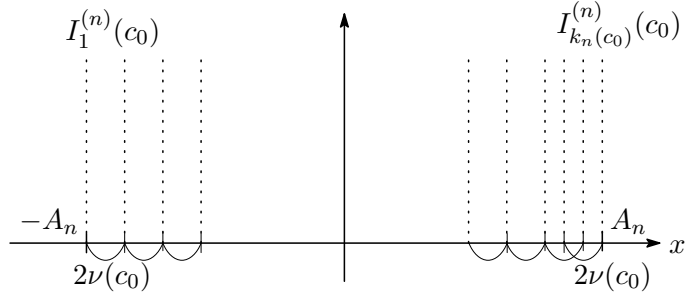}
  \caption{Division of $J_{0}^{(n)}$ by 
  short intervals of length $2\nu(c_{0})$.}
  \label{divided_J0n:tpic:fig}
 \end{center}
\end{figure}
Let $\epsilon > 0$ be arbitrarily fixed 
and let $J_{0}^{(n)} \equiv [-A_{n}, A_{n}]$.
We divide $J_{0}^{(n)}$ from $-A_{n}$ to $A_{n}$ by short intervals of 
length $2\nu(c_{0})$.  In right end of the intervals of $J_{0}^{(n)}$,
overlap of two short intervals of length $2\nu(c_{0})$ is allowed and the
right end of a short interval coincides with the right end of $J_0^{(n)}$.
See Figure \ref{divided_J0n:tpic:fig}.
Let $k_{n}({c_{0}})$ be the number of
short intervals and let
$I_{1}^{(n)}(c_{0}),\ldots,I_{k_{n}(c_{0})}^{(n)}(c_{0})$ be the divided short
intervals.
Then we have
\begin{eqnarray}
 k_{n}(c_{0})
  \leq 
   \frac{2A_{n}}{2\nu(c_{0})} + 1
   = 
   \frac{A_{n}}{\nu(c_{0})} + 1
   =
   \frac{A_{0} \cdot n^{\frac{2 + \zeta}{\beta - 1}}}{\nu(c_{0})} + 1
  \; .
  \label{bounding-number-of-short-intervals-c-zero:eq}
\end{eqnarray} 
Note that any interval in $J_{0}^{(n)}$ of length $2\nu(c_{0})$ is 
covered by at most $3$ small intervals 
from $\{I_{1}^{(n)}(c_{0}),\ldots,I_{k_{n}(c_{0})}^{(n)}(c_{0})\}$.
Now consider 
$J(\theta) = \bigcup_{m=1}^{K}
[\locparam_{m}-\nu({\scparam}_{m}),\locparam_{m}+\nu({\scparam}_{m}))$.
Since 
$[\locparam_{m}-\nu({\scparam}_{m}),\locparam_{m}+\nu({\scparam}_{m}))
\; , \; (m=1,\ldots,K)$ ,  
are intervals of length less
than or equal to  $2\nu(c_{0})$, 
$J(\theta)$ is covered by at most $3M$ short intervals.
Then 
\begin{eqnarray}
 \lefteqn{
  \sup_{\theta \in \Theta_{n,\scrk,s}'}
  \frac{1}{n} R_n(J(\theta))
  - 3M \cdot u_{0} \cdot 2\nu(c_{0})
  > \epsilon
  } & &
  \nonumber \\
 & \Rightarrow &
  \{x_{n,1} < - A_{n} \quad \rmor \quad x_{n,n} > A_{n}\}
  \nonumber \\
  & & \rmor 
  \nonumber \\
 & & 
  \{{1} \leq \exists k \leq {k_{n}(c_{0})} \; , \;  
  \frac{1}{n} R_n(I_{k}(c_{0}))
  - u_{0} \cdot 2\nu(c_{0})
  > \frac{\epsilon}{3M}
  \}
  \; . 
  \label{supR_InequalityUsingSmallIntervals:eq}
\end{eqnarray}
By lemma \ref{expanding-A:lem}, 
$\sum_n \prob(x_{n,1} < - A_{n} \  \rmor \  x_{n,n} > A_{n}) <
\infty$ and the first event on  the right-hand side of 
(\ref{supR_InequalityUsingSmallIntervals:eq}) can be ignored. 
We only need to consider the second event.
We will use the same logic in the proofs of  lemmas 
 \ref{99-2:goal:lem} and  \ref{99-1:goal:lem} below.
Then
\begin{eqnarray}
 & & 
 \prob
  \left(
   \sup_{\theta \in \Theta_{n,\scrk,s}'}
  \frac{1}{n} R_n(J(\theta))
  - 3M \cdot u_{0} \cdot 2\nu(c_{0})
  > \epsilon
  \right)
 \nonumber \\
 & & \hspace{6cm} \leq 
  \sum_{k=1}^{k_{n}(c_{0})}
  \prob
  \left(
   \frac{1}{n} R_n(I_{k}(c_{0}))
   - u_{0} \cdot 2\nu(c_{0})
   > \frac{\epsilon}{3M}
  \right)
  \; . 
  \nonumber 
\end{eqnarray}
Recall that, for any set $V \subset \real$, 
we denote by $P_0(V)$ the probability of
$V$ under the true density in (\ref{eq:def:P_0}) 
and denote by $R_{n}$ the number of observations
which belong to $V$ as in lemma~\ref{lem:fourterms}. 
Since 
\begin{eqnarray}
 P_0(I_{k}(c_{0})) \leq u_{0} \cdot 2\nu(c_{0}) , 
  \quad 
  (k=1,\ldots,k_{n}(\theta)), 
  \nonumber
\end{eqnarray} 
$R_n(V) \sim \bin(n , P_0(V))$ 
and from Okamoto's inequality (\cite{O1958}), 
we obtain \begin{eqnarray}
 \lefteqn{
  \prob
  \left(
   \frac{1}{n} R_n(I_{k}(c_{0}))
   - u_{0} \cdot 2\nu(c_{0})
   > \frac{\epsilon}{3M}
  \right)
  } & & \nonumber \\
 & \leq &
  \prob
  \left(
   \frac{1}{n} R_n(I_{k}(c_{0}))
   - P_0(I_{k}(c_{0}))
   > \frac{\epsilon}{3M}
  \right)
  \nonumber \\
 & \leq & 
  \exp{
  \left(
   -\frac{2n\epsilon^2}{9M^{2}}
  \right)} .
  \nonumber 
\end{eqnarray}
Therefore from (\ref{bounding-number-of-short-intervals-c-zero:eq}) 
\begin{eqnarray}
 \lefteqn{
 \prob
  \left(
   \sup_{\theta \in \Theta_{n,\scrk,s}'}
  \frac{1}{n} R_n(J(\theta))
  - 3M \cdot u_{0} \cdot 2\nu(c_{0})
  > \epsilon
  \right)
  } & & \nonumber \\
 & &  \hspace{4cm} 
  \le
  \left(
   \frac{A_{0} \cdot n^{\frac{2 + \zeta}{\beta - 1}}}{\nu(c_{0})} + 1
  \right)
  \cdot
  \exp{
  \left(
   -\frac{2n\epsilon^2}{9M^{2}}
  \right)}.
  \nonumber 
\end{eqnarray}
When we sum this over $n$, 
the resulting series on the right converges.
Hence by the Borel-Cantelli lemma, we have 
\begin{eqnarray}
 \prob
  \left(
    \sup_{\theta \in \Theta_{n,\scrk,s}'}
    \frac{1}{n} R_n(J(\theta))
    - 3M \cdot u_{0} \cdot 2\nu(c_{0})
   > \epsilon
   \quad i.o.
  \right) = 0.
  \nonumber 
\end{eqnarray}
Because $\epsilon > 0$ was arbitrary, 
we obtain 
\begin{eqnarray}
 \limsup_{n \rightarrow \infty}
 \sup_{\theta \in \Theta_{n,\scrk,s}'}
  \frac{1}{n} R_n(J(\theta))
  \leq
  3M \cdot u_{0} \cdot 2\nu(c_{0}), 
  \quad a.e. 
  \nonumber
\end{eqnarray}
\qed

\bigskip
By this lemma 
and 
(\ref{eq:condition:proof:mainthm})
we have
\begin{eqnarray}
  \limsup_{n \rightarrow \infty}
  \sup_{\theta \in \Theta_{n,\scrk,s}'}
  \frac{1}{n} R_n(J(\theta)) \cdot 
  \left(-\log{{\kappa_{0}}}\right)
  \leq 3M \cdot u_{0} \cdot 2\nu(c_{0}) \cdot
  \abslr{\log{{\kappa_{0}}}}
  < {\lambda_{0}}
  \quad a.e. 
 \label{boundingRJ:eq}
\end{eqnarray}
This bounds the third term on the right-hand side of 
(\ref{loglikelihood:lem:eq}) from above.

\subsubsection{Bonding the fourth term}

Finally 
we bound the fourth term on the right-hand side of
(\ref{loglikelihood:lem:eq}) from above. 
From lemma~\ref{lem:BoundMixtureByStepFunc} we have 
\begin{equation}
 \indicator_{J(\theta)}(x) \cdot 
 \sum_{m=1}^{M}
 \wtparam_{m} f_{m}(x; \locparam_{m}, \scparam_{m})
 \le 
 \sum_{t=1}^{T(\theta)}H(J_{t}(\theta)) \cdot
 \indicator_{J_{t}(\theta)}(x) 
 \quad , \quad 
 (x \in J(\theta))
 \label{bounded-by-uniform-density-2:main-thm:eq}
\end{equation}

\begin{figure}[htbp]
 \begin{center}
  \includegraphics{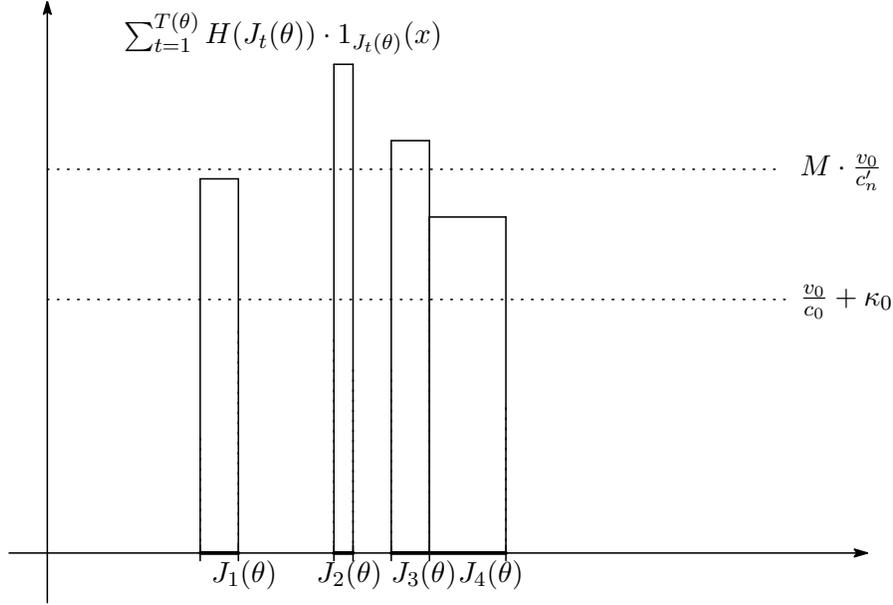}
  \caption{Example of classification of the intervals 
  by the height ($\tau_{n}(\theta) = \{1,4\} \; , \; \tau_{n}'(\theta) = \{2,3\} $).}
  \label{fig:CnCnd}
 \end{center}
\end{figure}

We now classify the intervals $J_{t}(\theta),\ t=1,\ldots, T(\theta),$
by the height $H(J_t(\theta))$.
Let 
\begin{eqnarray}
 c_{n}' \equiv c_{0} \cdot \exp{(-n^{1/4})}
  \label{eq:def:c_n_pr}
\end{eqnarray}
and define $\tau_n(\theta)$ and $\tau_n'(\theta)$ by 
\begin{eqnarray}
 \tau_n(\theta) \equiv
 \{
 t \in \{1,\ldots,T(\theta)\}
 \mid
 H(J_{t}(\theta)) \leq M v_0/c_{n}'
 \} 
\ , \ 
 \tau_n'(\theta) \equiv 
 \{1,\ldots,T(\theta)\}\backslash \tau_n(\theta)
 \ . 
 \label{eq:def:tau_n_and_tau_n_pr}
\end{eqnarray}
See Figure \ref{fig:CnCnd}.

Now suppose that the following inequality holds.
\begin{eqnarray}
 \lefteqn{
  \limsup_{n \rightarrow \infty}
  \sup_{\theta \in \Theta_{n,\scrk,s}'}
  \left[
   \sum_{t=1}^{T(\theta)}
   \frac{1}{n}
   R_{n}(J_{t}(\theta))
   \log{ H(J_{t}(\theta)) }
   \right.
  } & & \nonumber \\
 & & 
  - \left.
   3\left\{
     \sum_{t \in \tau_n(\theta)}
     u_{0} \cdot {\xi(H(J_{t}(\theta)))} \cdot
     \log{H(J_{t}(\theta))}
     +
     \sum_{t \in \tau_n'(\theta)}
     \frac{2}{n}\log{H(J_{t}(\theta))}
    \right\}
  \right]
  \leq 0 , 
  \quad a.e.
  \nonumber \\
  \label{2:goal:eq}
\end{eqnarray}
From 
(\ref{eq:condition:proof:mainthm}), 
and noting that $\log y/y^{\tilde{\beta}}$ is
decreasing in $y^{\tilde{\beta}} \ge e$, 
we have 
\begin{eqnarray}
 3 \cdot \sum_{t \in \tau_n(\theta)}
     {u_{0}} \cdot {\xi(H(J_{t}(\theta)))} \cdot 
     \log{H(J_{t}(\theta))}
   & \leq &
    3 \cdot 2M \cdot u_{0} \cdot 
    \xi(v_0/c_0)
    \cdot \log(v_0/c_{0}) 
   < 
  {\lambda_{0}} , 
  \nonumber \\
 3 \cdot \sum_{t \in \tau_n'(\theta)}
  \frac{2}{n}\log{H(J_{t}(\theta))}
  &\leq&
  3 \cdot 2{M} \cdot \frac{2}{n} \cdot 
  \log\frac{Mv_0}{c_n}
  \rightarrow 0 .
  \label{1:limit:ineq:eq}
\end{eqnarray}
Then from 
(\ref{bounded-by-uniform-density-2:main-thm:eq}), 
(\ref{2:goal:eq})
and
(\ref{1:limit:ineq:eq}), 
the fourth term on the right-hand side of 
(\ref{loglikelihood:lem:eq}) is bounded from above as
\begin{eqnarray}
\label{eq:bound-4nd-term:eq}
\limsup_{n\rightarrow\infty}\frac{1}{n} \sup_{\theta\in
  \Theta_{n,\scrk,s}'}
\sum_{x_i \in J(\theta)} \log
f(x_i; \theta) \le {\lambda_{0}} 
\quad a.e. 
\end{eqnarray}
Combining 
(\ref{loglikelihood:lem:eq}), 
(\ref{bounding-barf:eq}), 
(\ref{bounding-A-zero-term:eq}),
(\ref{boundingRJ:eq}) and
(\ref{eq:bound-4nd-term:eq}) we obtain 
\begin{eqnarray}
  \limsup_{n \rightarrow \infty}
  \sup_{\theta \in \Theta_{n,\scrk,s}'}
  \frac{1}{n}  \sum_{i=1}^{n} 
  \log{f(x_i; \theta)}
 & \leq &
E_{0}[\log{f(x;\theta_0)}]
   - {\lambda_{0}} ,
   \quad a.e.
   \nonumber 
\end{eqnarray}
and ~(\ref{goal:eq}) is satisfied.  Therefore it suffices to
prove (\ref{2:goal:eq}), which is a new goal of our proof.

We now consider further finite covering of $\Theta_{n,\scrk,s}'$.  For
any $T \; (1 \leq T \leq 2M)$ and $\tau \subset
\{1,\dots,T(\theta)\}$, define a subset of 
$\Theta_{n,\scrk,s}$ by
\begin{eqnarray}
 \Theta_{n,\scrk,s,T,\tau}'
  \equiv
  \{
  \theta \in \Theta_{n,\scrk,s}'
  \mid 
  T(\theta) = T 
  \; ,\; 
  \tau_{n}(\theta) = \tau 
  \}
 \; .
 \label{eq:def:Theta_scrk_s_T_tau}
\end{eqnarray}
Then (\ref{2:goal:eq}) is derived from the following two lemmas. 
\begin{lem}
 \begin{eqnarray}
  \limsup_{n \rightarrow \infty}
   \sup_{\theta \in \Theta_{n,\scrk,s,T,\tau}'}
   \left[
    \sum_{t \in \tau'}
    \frac{1}{n}
    R_{n}(J_{t}(\theta)) \cdot 
    \log{ H(J_{t}(\theta)) }
   -
   3 \sum_{t \in \tau'}
   \frac{2}{n}\log{H(J_{t}(\theta))}
  \right]
  \leq 0
  \quad a.e. \ , 
  \nonumber \\
  \nonumber 
 \end{eqnarray}
where $\tau' = \{1,\dots,T\}\backslash \tau$.
\label{99-2:goal:lem}
\end{lem}
\Proof
Let $\delta > 0$ be  any fixed positive real constant and  
let $\locparam_{t}'(\theta)$ denote the middle point of $J_{t}(\theta)$. 
Here, we consider the probability of the event that 
\begin{eqnarray}
 \sup_{\theta \in \Theta_{n,\scrk,s,T,\tau}'}
  \left[
   \sum_{t \in \tau'}
   \frac{1}{n}
   R_{n}(J_{t}(\theta)) \cdot 
   \log{ H(J_{t}(\theta)) }
   -
   3 \sum_{t \in \tau'}
   \frac{2}{n}\log{H(J_{t}(\theta))}
  \right]
  > 2M\delta .
  \label{99-2:RareEvent:eq}
\end{eqnarray}

Since 
$H(J_{t}(\theta)) > M v_0/c_{n}'$
holds for $t \in \tau'$, 
we obtain by lemma \ref{lem:BoundWidthOfInterval}
\begin{eqnarray}
 W(J_{t}(\theta))
 & \leq &
  v_{2} \cdot
  \left(
   \frac{c_{n}'}{M v_0}
  \right)
  ^{\tilde{\beta}}
  =
  v_{2} \cdot
  \left(
   \frac{c_{0}}{M v_{0}}
  \right)
  ^{\tilde{\beta}}
  \cdot \exp{(- {\tilde{\beta}} \cdot n^{1/4})}.
  \nonumber 
\end{eqnarray}
Let 
\begin{eqnarray}
 v_{3} & \equiv & 
  v_{2} \cdot
  \left(
   \frac{c_{0}}{M v_{0}}
  \right)
  ^{\tilde{\beta}}, 
  \nonumber \\
 w_{n} & \equiv & \frac{v_{3}}{2} \cdot \exp{(- {\tilde{\beta}} \cdot n^{1/4})}, 
  \label{eq:def:w_n}
  \\ 
  R_{n}[\mu, w] & \equiv & R_{n}([\mu - w, \mu + w]). 
   \label{eq:def:R_n-mu-w}
\end{eqnarray}
Noting that for $t \in \tau'$, the length of $J_t(\theta)$ 
is less than or equal to $2w_{n}$, 
the following relation holds. 
\begin{eqnarray}
 \lefteqn{
  \textrm{The event (\ref{99-2:RareEvent:eq}) occurs. }
  } & & 
  \nonumber \\
 & \Rightarrow &
 \sup_{\theta \in \Theta_{n,\scrk,s,T,\tau}'}
  \left[
   \sum_{t\in \tau'}
    \left(
     \frac{1}{n}
      R_n[\locparam_{t}'(\theta), w_{n}]
     - 3 \cdot \frac{2}{n}
    \right)
    \cdot 
   \log \frac{M v_0}{c_{n}}
  \right] 
 > 2M\delta
 \nonumber \\
 & \Rightarrow & 
  \exists \theta \in \Theta_{n,\scrk,s,T,\tau}' ,\; 
  \exists t \in \tau' 
  \quad \text{s.t.} \quad  
   \left(
    \frac{1}{n}
    R_n[\locparam_{t}'(\theta), w_{n}]
    - 3 \cdot \frac{2}{n}
   \right)
   \cdot 
  \log\frac{Mv_0}{c_{n}}
  > \delta 
  \nonumber \\
 & \Rightarrow &
  \exists \theta \in \Theta_{n,\scrk,s,T,\tau}' ,\; 
  \exists t \in \tau' 
  \quad \text{s.t.} \quad  
  R_n[\locparam_{t}'(\theta), w_{n}]
  \geq {6}
  \nonumber \\ 
 & \Rightarrow &
  \sup_{-\infty < {\locparam'} < \infty} R_n[{\locparam'}, w_{n}] \geq 6
   \; .
  \label{99-2:rarer-event:eq}
\end{eqnarray}
Below, we consider the probability of the event 
that (\ref{99-2:rarer-event:eq}) occurs. 
We divide $J_{0}^{(n)} = [-A_{n}, A_{n}]$
 from $-A_{n}$ to $A_{n}$ 
by short intervals of length $2w_{n}$ as 
in the proof of lemma \ref{boundedRJ:lem}.
Let $k({w_{n}})$ be the number of
short intervals and let
$I_{1}(w_{n}),\ldots,I_{k(w_{n})}(w_{n})$ be the divided short
intervals. 
Then we have
\begin{eqnarray}
 k(c_{n}') \leq 
   \frac{2A_{n}}{2w_{n}} + 1
   =
   \frac{A_{0} \cdot n^{\frac{2 + \zeta}{\beta - 1}}}{\nu(c_{0})} + 1
   \; .
  \label{99-2:bound-k_c_n:eq}
\end{eqnarray}  
Since any interval in $J_{0}$ of length $2w_{n}$ is 
covered by at most $3$ small intervals 
from $I_{1}(w_{n}),\ldots,I_{k(w_{n})}(w_{n})$ 
and from lemma \ref{expanding-A:lem},  
\begin{eqnarray}
 \sup_{-\infty < {\locparam'} < \infty} R_n[{\locparam'}, w_{n}] \geq 6
  & \Rightarrow &
  1 \leq \exists k \leq k(w_{n}) \; ,\; 
   R_n(I_{k}(w_{n})) \geq 2
  \; . 
  \nonumber 
\end{eqnarray}
Note that $R_n(I_{k}(w_{n}))
\sim \bin(n , P_0(I_{k}(w_{n})))$
and $P_0(I_{k}(w_{n})) \le 2w_{n} u_{0}$.
Therefore from (\ref{99-2:bound-k_c_n:eq}) 
we have 
\begin{eqnarray}
 \lefteqn{
  \sum_{k=1}^{k(w_{n})}
  \prob
  \left(
   R_n(I_{k}(w_{n})) \geq 2
  \right)
  \leq
  \left(
   \frac{A_{n}}{w_{n}} + 1
  \right) \cdot 
  \left\{
  \max_{1\le k \le k(w_{n})}
  \prob(R_n(I_{k}(w_{n})) \geq 2)
  \right\}
  } & & 
   \nonumber \\
 &  & 
  \leq
  \left(
   \frac{A_{n}}{w_{n}} + 1
  \right)
  \sum_{k=2}^{n}
  \begin{pmatrix} n \\ k \end{pmatrix}
  (2w_{n} u_{0})^{k}(1 - 2w_{n} u_{0})^{n-k}
   \nonumber \\
 &  &
  \leq
  \left(
   \frac{A_{n}}{w_{n}} + 1
  \right)
  \sum_{k=2}^{n}\frac{n^{k}}{k!}(2w_{n} u_{0})^{k}
  \leq
   \left(
   \frac{A_{n}}{w_{n}} + 1
  \right)
  (2nw_{n} u_{0})^{2}
  \sum_{k=0}^{n}\frac{1}{k!}(2nw_{n} u_{0})^{k}
  \nonumber \\
 &  &
  \leq
    \left(
   \frac{A_{n}}{w_{n}} + 1
  \right)
  (2nw_{n} u_{0})^{2}
  \exp{(2nw_{n} u_{0})}
  \; .
  \nonumber 
\end{eqnarray}
When we sum this over $n$, 
resulting series on the right converges.
Hence by the Borel-Cantelli lemma and the fact that $\delta > 0$ was arbitrary, 
we obtain 
\begin{eqnarray}
  \limsup_{n \rightarrow \infty}
   \sup_{\theta \in \Theta_{n,\scrk,s,T,\tau}'}
   \left[
    \sum_{t \in \tau'}
    \frac{1}{n}
    R_{n}(J_{t}(\theta)) \cdot 
    \log{ H(J_{t}(\theta)) }
    -
    3 \sum_{t \in \tau'}
    \frac{2}{n}\log{H(J_{t}(\theta))}
  \right]
   \leq 0
   \quad a.e.
   \nonumber 
 \end{eqnarray}
\qed 

\begin{lem}
 \begin{eqnarray}
  \lefteqn{
   \limsup_{n \rightarrow \infty}
   \sup_{\theta \in \Theta_{n,\scrk,s,T,\tau}'}
   \left[
    \sum_{t \in \tau}
    \frac{1}{n}
    R_{n}(J_{t}(\theta)) \cdot 
    \log{ H(J_{t}(\theta)) }
   \right.
   }& & 
   \nonumber \\
  & & \hspace{4cm}\left.
   -
   3 \sum_{t \in \tau}
   u_{0} \cdot \xi(H(J_{t}(\theta)))
   \cdot \log{H(J_{t}(\theta))}
  \right]
  \leq 0
  \quad a.e.
  \nonumber
 \end{eqnarray}
 \label{99-1:goal:lem}
\end{lem}
\Proof 
Let $\delta > 0$ be any fixed positive real constant 
and let 
\begin{eqnarray}
 h_{n}
  \equiv
  \frac{\delta}{12}
  \left\{
   u_{0} \cdot 
   \log{
     \frac{M v_0}{c_{n}'}
   }
  \right\}^{-1}
 \; .
 \label{99-1:h_n:def:eq}
\end{eqnarray}
Since 
$v_0/c_0 \leq H(J_{t}(\theta)) \leq M v_0/c_{n}'$,
we have 
$\xi(M v_0/c_{n})\leq \xi(H(J_{t}(\theta))) \leq \xi(v_0/c_{0})$.
We divide the interval 
$[\xi(M v_0/c_{n}'), \xi(v_0/c_{0})]$ 
from 
$\xi(c_0/v_0)$ to $\xi(M v_0/c_{n}')$ 
by short intervals of length $h_{n}$. 
In the left end 
$\xi(M v_0/c_{n}')$ 
of the interval 
$[\xi(M v_0/c_{n}'), \xi(v_0/c_{0})]$,
overlap of two short intervals of length $h_{n}$ is allowed and the
left  end of a short interval is equal to 
$\xi(M v_0/c_{n}')$.
Let $l_{n}$ be the number of
short intervals of length $h_{n}$ and define $w_{l}^{(n)}$ by
\begin{eqnarray}
 2w_{l}^{(n)} \equiv 
 \begin{cases}
   \xi(v_0/c_{0}) - (l - 1)h_{n}, & 1 \leq l \leq l_{n},\\
    \xi(M v_0/c_{n}'), & l = l_{n} + 1.
 \end{cases}
 \label{eq:def:w_l-n}
\end{eqnarray}
Then we have 
\begin{eqnarray}
 l_{n} \leq 
   \frac{\xi(v_0/c_{0})}{h_{n}} + 1
  \; . 
  \label{99-1:l_n:ineq:eq}
\end{eqnarray}
Let 
\begin{eqnarray}
 \psi(y) \equiv \xi^{-1}(y) = \left(\frac{v_{2}}{y}\right)^{1/\tilde{\beta}} \ , 
  \qquad (y > 0) \ ,
  \nonumber 
\end{eqnarray}
where $\xi^{-1}(\cdot)$ is the inverse function of $\xi(\cdot)$.
Next we consider the probability of the event that 
\begin{eqnarray}
 \lefteqn{
  \sup_{\theta \in \Theta_{n,\scrk,s,T,\tau}'}
  \left[
   \sum_{t \in \tau}
   \frac{1}{n}
   R_{n}(J_{t}(\theta)) \cdot 
   \log{ H(J_{t}(\theta)) }
  \right.
   }& & 
   \nonumber \\
 & & \hspace{4cm}
  \left.
   -
   3 \sum_{t \in \tau}
   u_{0} \cdot \xi(H(J_{t}(\theta)))
   \cdot \log{H(J_{t}(\theta))}
  \right]
  > 2M\delta .
  \label{99-1:rare-event:eq}
\end{eqnarray}
For this event the following relation holds.
\begin{eqnarray}
 \lefteqn{
  \textrm{The event (\ref{99-1:rare-event:eq}) occurs. }
  } & & 
  \nonumber \\
 & \Rightarrow &
 \exists \theta \in \Theta_{n,\scrk,s,T,\tau}' 
 , \; 
  \forall t \in \tau ,\; 
  1 \leq \exists l(t) \leq l_{n} 
  \quad {\rm s.t.} \; 
  \nonumber \\
 & &
 \psi(2w_{l(t)}^{(n)}) \leq H(J_{t}(\theta))
 \leq \psi(2w_{l(t)+1}^{(n)})
 \quad {\rm and} \;  
 \nonumber \\
 & & 
 \sum_{t \in \tau}
   \left(
    \frac{1}{n}
    R_n[\locparam_{t}'(\theta), w_{l(t)}^{(n)}]
   -
   3u_{0} \cdot 2w_{l(t)+1}^{(n)}
  \right)
  \cdot 
 \log\psi(2w_{l(t)+1}^{(n)})
 > 2M\delta
 \nonumber 
 \\
 & \Rightarrow &
  \exists \theta \in \Theta_{n,\scrk,s,T,\tau}' 
  \; , \; 
  \exists t \in \tau 
  , \; 
  1 \leq 
  \exists l(t)
  \leq l_{n}
  \quad {\rm s.t.} \; 
  \nonumber \\
 & &
  \psi(2w_{l(t)}^{(n)}) \leq H(J_{t}(\theta)) \leq \psi(2w_{l(t)+1}^{(n)})
  \quad {\rm and} \;  
  \nonumber \\
 & &
   \left(
    \frac{1}{n}
    R_n[\locparam_{t}'(\theta), w_{l(t)}^{(n)}]
      -
    3u_{0} \cdot 2w_{l(t)+1}^{(n)}
   \right)
   \cdot 
  \log\psi(2w_{l(t)+1}^{(n)})
  > \delta
  \nonumber \\
   & \Rightarrow &
  1 \leq 
  \exists l
  \leq l_{n}
  \quad {\rm s.t.} \; 
  \nonumber \\
 & &
   \sup_{-\infty < {\locparam'} < \infty}
   \left(
    \frac{1}{n}
    R_n[{\locparam'}, w_{l}^{(n)}]
      -
    3u_{0} \cdot 2w_{l+1}^{(n)}
   \right)
   \cdot 
  \log\psi(2w_{l+1}^{(n)})
  > \delta
  \nonumber \\
 & \Rightarrow &
  1 \leq \exists l \leq l_{n}
  \quad {\rm s.t.} \; 
  \nonumber \\
 & &
  \sup_{-\infty < {\locparam'} < \infty}
  \left\{
  \left(
   \frac{1}{n} R_n[{\locparam'}, w_{l}^{(n)}]
   -
   3u_{0} \cdot 2w_{l}^{(n)}
   \right) \cdot 
  \log\psi(2w_{l+1}^{(n)})
  \right.
  \nonumber \\
 & & \qquad \qquad \qquad \qquad \quad 
  +
  \left.
   3u_{0} \cdot (2w_{l}^{(n)} - 2w_{l+1}^{(n)})
   \cdot \log \psi(2w_{l+1}^{(n)})
  \right\}
  > \delta
 \label{99-1:rarer-event:eq}
\end{eqnarray}
Then from (\ref{99-1:h_n:def:eq}) and lemma \ref{expanding-A:lem} 
the following relation holds.
\begin{eqnarray}
 \lefteqn{
  \textrm{The event (\ref{99-1:rarer-event:eq}) occurs.}
  } & &
 \nonumber \\
 & \Rightarrow &
  {1}\leq \exists l \leq {l_{n}} 
  \quad {\rm s.t.} \; 
  \nonumber \\
 & & 
   \sup_{-\infty < {\locparam'} < \infty}
   \frac{1}{n}
   \left(
    R_n[{\locparam'}, w_{l}^{(n)}]
    -
    3u_{0} \cdot 2w_{l}^{(n)}
   \right) \cdot 
   \log \psi(2w_{l+1}^{(n)})
   > \frac{\delta}{2}
  \nonumber \\
  \nonumber \\
 & \Rightarrow &
  {1}\leq \exists l \leq {l_{n}} 
  \quad {\rm s.t.} \; 
  \nonumber \\ 
 & & 
  \sup_{-A_{n} < {\locparam'} < A_{n}}
  \frac{1}{n}
   \left(
    R_n[{\locparam'}, w_{l}^{(n)}]
    -
    3u_{0} \cdot 2w_{l}^{(n)}
   \right) \cdot 
   \log \psi(2w_{l+1}^{(n)})
   > \frac{\delta}{2} 
  \label{99-1:rarest-event:eq}
\end{eqnarray}
Below, we consider the probability of the event 
that 
(\ref{99-1:rarest-event:eq}) occurs.
We divide $J_{0}^{(n)}$ from $-A_{n}$ to $A_{n}$ 
by short intervals of length $2w_{l}^{(n)}$ as 
in the proof of lemma \ref{boundedRJ:lem}. 
Let $k({w_{l}^{(n)}})$ be the number of
short intervals and let
$I_{1}(w_{l}^{(n)}),\ldots,I_{k(w_{l}^{(n)})}(w_{l}^{(n)})$ be the divided short
intervals. 
Then we have
\begin{eqnarray}
 k(w_l^{(n)}) \leq 
   \frac{2A_{n}}{2w_{l}^{(n)}} + 1
  \; .
  \label{99-1:bound-k_c_n:eq}
\end{eqnarray} 
Since any interval in $J_{0}$ of length $2\scparam_l^{(n)}$ is 
covered by at most $3$ small intervals 
from $\{I_{1}(w_{l}^{(n)}),\ldots,I_{k(w_{l}^{(n)})}(w_{l}^{(n)})\}$,
we have 
\begin{eqnarray}
 \lefteqn{
  \sup_{-A_{n} \leq {\locparam'} \leq A_{n}}
  \left(
   \frac{1}{n}
   R_n[{\locparam'}, w_{l}^{(n)}]
   -
   3u_{0} \cdot 2w_{l}^{(n)}
  \right)
  > 
  \frac{\delta}{2}
   \left(
    \log\psi(2w_{l(t)+1}^{(n)})
   \right)^{-1}
   } & & \nonumber \\
 & \Rightarrow &
  \max_{k = 1,\ldots,k(w_{l}^{(n)})}
  \left(
   \frac{1}{n}
   R_n(I_k(w_{l}^{(n)}))
   -
   u_{0} \cdot 2w_{l}^{(n)}
  \right)
  > 
  \frac{1}{3} \cdot \frac{\delta}{2}
  \left(
   \log\psi(2w_{l(t)+1}^{(n)})
  \right)^{-1}. 
  \label{99-1:short-interval-event:eq}
\end{eqnarray}
Note that $R_n(I_k(w_{l}^{(n)})) \sim \bin(n, P_0(I_k(w_{l}^{(n)})))$ 
and $P_0(I_k(w_{l}^{(n)})) \leq u_{0} \cdot 2w_{l}^{(n)}$.
Therefore from (\ref{99-1:bound-k_c_n:eq}) and Okamoto's inequality (\cite{O1958}) we have  
\begin{eqnarray}
 \lefteqn{
  \prob
  \left(
   \max_{k = 1,\ldots,k(w_{l}^{(n)})}
   \frac{1}{n}
   \left( 
    R_n(I_k(w_{l}^{(n)}))
    -
    u_{0} \cdot 2w_{l}^{(n)}
   \right)
   > 
   \frac{1}{3} \cdot \frac{\delta}{2}
   \left\{
    \log{\psi(2w_{l(t)+1}^{(n)})}
   \right\}^{-1}
  \right)
  } & & \nonumber \\
  & & \hspace{1.5cm} 
   \leq 
   \left(
    \frac{2A_{n}}{2w_{l}^{(n)}} + M
   \right) \cdot 
   \exp{
   \left[
    -2n \cdot \frac{\delta^{2}}{36}
    \left\{
     \log{\psi(2w_{l(t)+1}^{(n)})}
    \right\}^{-2}
   \right]
   }
   \nonumber \\
 & & \hspace{1.5cm} 
   \leq 
   \left(
    \frac{A_{n}}{\xi(M v_0/c_{n}')} + M
   \right)
   \cdot
   \exp{
   \left[
    -2n \cdot \frac{\delta^{2}}{36}
    \left\{
     \log(M v_0/c_{n}')
    \right\}^{-2}
   \right]
   }
   \label{99-1:c_n_dash:prob:eq}
   . 
\end{eqnarray}
{}From (\ref{99-1:l_n:ineq:eq}), 
(\ref{99-1:rarer-event:eq}), 
(\ref{99-1:rarest-event:eq}),
(\ref{99-1:short-interval-event:eq}),
and (\ref{99-1:c_n_dash:prob:eq})
we obtain 
\begin{eqnarray}
 & & 
  \sum_{l=1}^{l_{n}}
  \prob
  \left(
   \sup_{-A_{n} < {\locparam'} < A_{n}}
   \frac{1}{n}
   \left(
    R_n[{\locparam'}, w_{l}^{(n)}]
    -
    3u_{0} \cdot 2w_{l}^{(n)}
   \right) \cdot 
   \log \psi(2w_{l+1}^{(n)})
   > \frac{\delta}{2}
  \right)
  \nonumber 
\\&&  \quad 
\le 
\left(\frac{\xi(v_0/c_{0})}{h_{n}} + 1\right)
\cdot
\left(
    \frac{A_{n}}{\xi(M v_0/c_{n}')} + M
   \right)
   \cdot
   \exp{
   \left[
    -2n \cdot \frac{\delta^{2}}{36}
    \left\{
     \log{\left(M v_0/c_{n}'\right)}
    \right\}^{-2}
   \right]}
  \nonumber
\end{eqnarray}
When we sum this over $n$, 
the resulting series on the right converges.
Hence by the Borel-Cantelli lemma and the fact 
that $\delta > 0$ is arbitrary, we have
\begin{eqnarray}
 \lefteqn{
  \limsup_{n \rightarrow \infty}
  \sup_{\theta \in \Theta_{n,\scrk,s,T,\tau}'}
  \left[
   \sum_{t \in \tau}
   \frac{1}{n}
   R_{n}(J_{t}(\theta)) \cdot 
   \log{ H(J_{t}(\theta)) }
  \right.
  }& & 
  \nonumber \\
 & & \hspace{3cm}
  \left.
   -
   3 \sum_{t \in \tau}
   {u_{0}} \cdot \xi({H(J_{t}(\theta))})
   \cdot \log{H(J_{t}(\theta))}
  \right]
  \leq 0
  \quad a.e.
  \nonumber
\end{eqnarray}
\qed

This completes the proof of theorem \ref{main-thm:thm}.

\section{Discussions}
\label{sec:discussion}

In this paper we consider the strong consistency 
of MLE for mixtures of 
location-scale distributions.
We treat the case that the scale parameters of 
the component distributions are restricted from below 
by $c_n = \exp(-n^d)$, $0 < d < 1$, and 
give the regularity conditions for the strong consistency 
of MLE.

As in the case of the uniform mixture in \cite{TT2003-20}, it is readily verified that 
if $c_n$ decreases to zero faster than $\exp(-n)$, 
then the consistency of MLE  fails.
Therefore the rate of $c_{n} = \exp(-n^d)$, $0 < d < 1$, 
obtained in this paper 
is almost the lower bound of the order of $c_n$ which maintains
the strong consistency.

Although we treat the univariate case in this paper, it is clear that 
the result obtained in this paper 
can be extended to the multivariate case 
under the condition that 
components 
are bounded and their tails decrease to zero fast enough 
if the minimum singular values of the scale matrices 
of the components are restricted from below by $c_{n}$.

Finally let us consider some 
sufficient conditions for the regularity conditions. 
For $\theta_{m} \in \Omega_{m}$ and any positive real number $\rho$, let 
\begin{eqnarray}
f_{m}(x;\theta_{m},\rho) & \equiv & \sup_{\dist(\theta_{m}',\theta_{m}) \leq \rho}f_m(x;\theta_{m}'). 
\nonumber
\end{eqnarray}
Let $\Gamma$ be any compact
subset of $\Omega_{m}$.
Consider the following two conditions.
\begin{assumption}
 For each $\theta_{m} \in \Gamma$ and sufficiently small $\rho$, 
 $f_{m}(x;\theta_{m}, \rho)$ is measurable.
 \label{assumption:5}
\end{assumption}
\begin{assumption}
 For each $\theta_{m} \in \Gamma$,
 if
 $\lim_{\seqnum \rightarrow \infty}\theta_{m}^{(\seqnum)} = \theta_{m}$, 
 then 
 $\lim_{\seqnum \rightarrow \infty}
 f_{m}(x;\theta_{m}^{(\seqnum)}) = f_{m}(x;\theta_{m})$
 for all $x$.
 \label{assumption:6}
\end{assumption}
If assumptions~\ref{assumption:5} and \ref{assumption:6} hold,
then it is easily verified that 
assumptions~\ref{assumption:2} and \ref{assumption:3}
hold.
Thus assumptions~\ref{assumption:1}, \ref{assumption:4}, \ref{assumption:5} and \ref{assumption:6} are  
sufficient conditions for regularity conditions 
and assumptions~\ref{assumption:5} and \ref{assumption:6} are checked more easily.
For example, finite mixture density which consists of 
normal density, $t$-density and uniform density on an open interval 
satisfies assumptions~\ref{assumption:1}, \ref{assumption:4}, \ref{assumption:5} and \ref{assumption:6}. 

\bibliographystyle{econometrica}
\bibliography{mixture}

\clearpage 
{ 
\thispagestyle{empty}
\begin{table}
\caption{List of notations \quad  ($\theta \in \Theta;\; \scrk \subset \{1,\dots,M\};\; V \subset \real ;\; y ,\rho , \mu, w \in \real$)}
\medskip
\begin{tabular}{|c|c|} \hline 
 \bf{Notation} & \bf{Definition or description} \\ \hline 
 $M$ & Number of components \\ \hline 
 $\theta_{\scrk}$ & Subvector of $\theta \in \Theta$ 
consisting of the components in $\scrk$ \\ \hline
 $\bar{\Theta}_{\scrk}$ & 
 $\Theta_{\scrk} \equiv \{\theta_{\scrk} \mid \theta \in \Theta\}$
 ; Parameter space of $\bar{\theta}_{\scrk}$ ; See (\ref{eq:def:barTheta_scrk}) \\ \hline
 $f_{\scrk}(x;\theta_{\scrk})$ & 
 $f_{\scrk}(x;\theta_{\scrk}) \equiv \sum_{k \in \scrk} \wtparam_{k} f_k(x;{\locparam}_k, \scparam_{k})$ ; See (\ref{eq:def:f_scrk})\\ \hline
 $f_{\scrk}(x;\theta_{\scrk}, \rho)$ & 
 $ f_{\scrk}(x;\theta_{\scrk},\rho) 
 \equiv 
 \sup_{\dist(\theta_{\scrk}',\theta_{\scrk}) \leq \rho}f_{\scrk}(x;\theta_{\scrk}')$ ; See (\ref{eq:def:f_scrk_rho})\\ \hline
 $\scrg_{\scrk}$ & 
 $\scrg_{\scrk} 
   \equiv 
  \{
  f_{\scrk}(x;\theta_{\scrk}) 
  \mid 
  \theta_{\scrk} \in \bar{\Theta}_{\scrk}
  \} $ ; See (\ref{eq:def:scrg_scrk})
 \\ \hline
 $\scrg_{K}$ & 
 $\scrg_{K} \equiv \bigcup_{\abslr{\scrk}\leq K}\scrg_{\scrk}$ ; See (\ref{eq:def:scrg_K})
 \\ \hline
 $v_{0}, v_{1}, \beta$ & 
 $ 
 f_{m}(x;\locparam_{m}=0, \scparam_{m}=1) 
 \leq 
 \min \{v_{0} \;, \; v_{1} \cdot \abs{x}^{-\beta}\}
 $ ; See Assumption~\ref{assumption:1}
 \\ \hline
 $c_{0}, c_{n}, d, \Theta_{n}$ & 
 $c_{n} = c_{0} \cdot \exp(-n^{d})$ ; 
 See theorem~\ref{main-thm:thm} \\ \hline
 $B$ & 
 $B \equiv {v_{0}}/{\kappa_{0}}$ ; See (\ref{eq:def:B})
 \\ \hline
 $\tilde{\beta}, \nu(y)$ & 
 $ 
 \tilde{\beta}
  \equiv 
  \frac{\beta-1}{\beta}, \; 
 \nu(y) \equiv 
  \left(
   \frac{v_{1}}{\kappa_{0}}
  \right)^{\frac{1}{\beta}}
  y^{\tilde{\beta}}
 $ ; See (\ref{eq:def:beta_and_nu_and_v})
 \\ \hline
 $J(\theta)$ & 
 $J(\theta)
  \equiv
 \bigcup_{m \in \scrk_{{\scparam} \leq c_{0}}}
 [\locparam_{m} - \nu(\scparam_{m}), \locparam_{m} + \nu(\scparam_{m}))$
 ; See (\ref{eq:def:Jtheta}) \\ \hline
 $J_{t}(\theta)$ & Interval of step function; See lemma~\ref{lem:BoundMixtureByStepFunc} \\ \hline
 $H(J_{t}(\theta))$ & Height of step function in $J_{t}(\theta)$ ; See lemma~\ref{lem:BoundMixtureByStepFunc} \\ \hline
 $W(J_{t}(\theta))$ & Width of $J_{t}(\theta)$; See lemma~\ref{lem:BoundMixtureByStepFunc} \\ \hline
 $T(\theta)$ & Number of steps ; See lemma~\ref{lem:BoundMixtureByStepFunc} \\ \hline
 $v_{2}, \xi(y)$ & 
 $ 
 v_{2}
   \equiv 
   2\left(
     \frac{v_{1}}{\kappa_{0}}
    \right)^{\frac{1}{\beta}}
   \left(
    {v_{0}{(M + 1)}}
   \right)^{\tilde{\beta}}, \; 
  \xi(y) \equiv v_{2}\cdot \left(\frac{1}{y}\right)^{\tilde{\beta}} 
 $ ; See (\ref{eq:def:v_2_and_xi})
 \\ \hline
 $u_{0}, u_{1}$ & 
 $ 
 f(x;\theta_{0})
  \leq
  \min{\{u_{0}, \; u_{1}\cdot \abs{x}^{-\beta}\}}
 $ ; See lemma~\ref{lem:BoundComponentByStepFunc} 
 \\ \hline
 $x_{n,1}, x_{n,n}$ & 
 $x_{n,1}  \equiv  \min{\{x_{1}, \ldots, x_{n}\}}, \;
 x_{n,n}  \equiv  \max{\{x_{1}, \ldots, x_{n}\}}$
 ; See lemma~\ref{expanding-A:lem}
 \\ \hline
 $A_{0}, \zeta, A_{n}$ & 
 $A_{n} \equiv A_{0} \cdot n^{\frac{2 + \zeta}{\beta - 1}}$ ; See lemma~\ref{expanding-A:lem} \\ \hline
 $A^{(j)}$ & 
 $ 
 \min_{m \in \scrk_{\abs{\locparam}\uparrow \infty}}
  \left\{
   \min{\{\abs{{\locparam}_{m}^{(\seqnum)} + \nu({\scparam}_{m}^{(\seqnum)})} , \abs{{\locparam}_{m}^{(\seqnum)} - \nu({\scparam}_{m}^{(\seqnum)})}\}}
  \right\} 
 $
 ; See (\ref{eq:def:A_hat_seqnum})
 \\ \hline
 $\scra_{0}$ & $\scra_{0} \equiv (-\infty, -A_{0}] \cup [A_{0},\infty)$ ; See (\ref{eq:def:scra_0})\\ \hline
 $\scrk_{{\scparam}\downarrow 0}, \scrk_{{\scparam}\uparrow \infty}, \scrk_{\abs{\locparam}\uparrow \infty}$ & 
 Disjoint subset of $\scrl$; See lemma~\ref{lem:proof-theor-refk} \\ \hline
 $\scrk_{{\scparam} \leq c_{0}}, \scrk_{{\scparam}\geq B}, \scrk_{\abs{\locparam} \geq A_{0}}$ & 
 Disjoint subset of $\{1,\dots,M\}$; See (\ref{eq:def:Theta_n_scrk})
 \\ \hline
 $\scrk_{R}$ & 
 \begin{tabular}{c} 
  $\scrkrest \equiv \scrl 
 \backslash \{\scrk_{{\scparam}\downarrow 0} \cup \scrk_{{\scparam}\uparrow \infty} 
 \cup \scrk_{\abs{\locparam}\uparrow \infty}\}$ in subsection~\ref{sec:proof-kappa-lambda-thm}
  \\ 
 $\scrkrest \equiv \{1, \dots, M\} 
 \backslash 
 \{\scrk_{{\scparam} \leq c_{0}} \cup \scrk_{{\scparam}\geq B} \cup \scrk_{\abs{\locparam} \geq A_{0}}\}$ in subsection~\ref{sec:notat-some-lemm} 
 \end{tabular} \\ \hline 
 $P_{0}(V)$ & 
 $P_0(V) \equiv \int_{V} f(x;\theta_{0}) \rmd x$
 ; See (\ref{eq:def:P_0})
 \\ \hline
 $R_{n}(V)$ & Number of observations which belong to a set $V$ \\ \hline
 $\scrb(\theta,\rho(\theta))$ & Open ball with center $\theta$ and radius $\rho(\theta)$ \\ \hline
 $c_{n}'$ & $c_{n}' \equiv c_{0} \cdot \exp{(-n^{1/4})}$ ; See (\ref{eq:def:c_n_pr}) \\ \hline
 $\tau_{n}(\theta), \tau_{n}(\theta)'$ & See (\ref{eq:def:tau_n_and_tau_n_pr}) \\ \hline
 $\Theta_{n}'$ & 
 $\Theta_{n}'
 \equiv
 \{
 \theta \in \Theta
 \mid
 \exists m \; {\rm s.t.} \; 
 c_{n} \leq \scparam_{m} \leq c_{0}
 \; \mathrm{or} \; 
 \abs{\locparam_{m}} > A_{0}
 \}$ \\ \hline
 $\Theta_{n, \scrk}'$ & See (\ref{eq:def:Theta_n_scrk})\\ \hline
 $\Theta_{n, \scrk, s}'$ & 
 $ 
 \Theta_{n,\scrk ,s}'
 \equiv 
 \{
 \theta \in \Theta_{n,\scrk}'
 \mid 
 \theta_{\scrkrest} \in \scrb(\theta_{\scrkrest}^{(s)},\rho(\theta_{\scrkrest}^{(s)}))
 $; See (\ref{eq:def:Theta_scrk_n_s})
 \\ \hline
 $\Theta_{n, \scrk, s, T, t}'$ & 
 $
  \Theta_{n,\scrk,s,T,\tau}'
  \equiv
  \{
  \theta \in \Theta_{n,\scrk,s}'
  \mid 
  T(\theta) = T 
  \; ,\; 
  \tau_{n}(\theta) = \tau 
  \}
 $ ; See (\ref{eq:def:Theta_scrk_s_T_tau}) \\ \hline
 $\locparam_{t}'(\theta)$ & 
 $\locparam_{t}'(\theta)$ denote the middle point of $J_{t}(\theta)$
 \\ \hline
 $w_{n}$ & See (\ref{eq:def:w_n})
  \\ \hline
 $w_{n}^{(l)}$ & See (\ref{eq:def:w_l-n})
  \\ \hline
 $R_{n}[\mu, w]$ & 
 $R_{n}[\mu, w] \equiv  R_{n}([\mu - w, \mu + w])$; See (\ref{eq:def:R_n-mu-w})
 \\ \hline
\end{tabular} 
\end{table}
}
\end{document}